\documentclass{article} 

\usepackage{latexsym}
\usepackage{amssymb}
\usepackage{amsfonts}
\usepackage{amsmath}
\usepackage{indentfirst}

\newcommand{\ff}{\varphi}
\renewcommand{\ll}{\lambda}

\renewcommand{\AA}{\mathcal{A}} 

\newcommand{\CC}{\mathcal{C}} 
\newcommand{\DD}{\mathcal{D}} 
\newcommand{\MM}{\mathcal{M}} 
\renewcommand{\SS}{\mathcal{S}} 
\newcommand{\YY}{\mathcal{Y}} 
\newcommand{\LL}{\mathcal{L}} 
\newcommand{\HH}{\mathcal{H}} 

\newcommand{\JJ}{\mathcal{J}} 

\newcommand{\NN}{\mathbb{N}} 
\newcommand{\ZZ}{\mathbb{Z}} 
\newcommand{\RR}{\mathbb{R}} 

\newcommand{\Spec}{\mathbf{Spec}} 
\newcommand{\Int}{\mathbf{Int}} 

\newcommand{\map}[3]{ { #1 } : { #2 } \rightarrow { #3 } }

\newtheorem{theorem}{Theorem}[section]
\newtheorem{proposition}[theorem]{Proposition}
\newtheorem{lemma}[theorem]{Lemma}
\newtheorem{corollary}[theorem]{Corollary}
\newtheorem{conjecture}[theorem]{Conjecture}
\newcommand{\eproof}{ \hspace*{ \fill } $ \Box $ \vspace{1mm} }
\newenvironment{proof}{ {\em Proof}. }{ \eproof }

\newenvironment{remark}{\medskip\textbf{Remark.}}{}

\title{Lattices of Paths: \\ Representation Theory and Valuations}
\author{Luca Ferrari\thanks{Dipartimento di Sistemi e Informatica, viale Morgagni 65, 50134 Firenze, Italy
{\tt ferrari@dsi.unifi.it}}
\and
Emanuele Munarini\thanks{Politecnico di Milano, Dipartimento di Matematica,
Piazza Leonardo da Vinci 32, 20133 Milano, Italy {\tt emanuele.munarini@polimi.it}}}
\date{}

\begin{document}

\maketitle

\begin{abstract}
 We study some distributive lattices arising in the combinatorics of lattice paths.
 In particular, for the Dyck, Motzkin and Schr\"oder lattices
 we describe the spectrum and we determine explicitly the Euler characteristic
 in terms of natural parameters of lattice paths.
\end{abstract}

\bigskip


\bigskip


\section{Introduction}

The set of all lattice paths of equal length with steps of a prescribed kind
(usually starting from the origin and ending on the $x$-axis with respect to a fixed Cartesian coordinate system)
can always be ordered by containment.
More precisely, we say that $\; \gamma_1 \leq \gamma_2 \;$ when $\; \gamma_1 \;$ lies weakly below $\; \gamma_2 \,$.
In several cases, the resulting poset has the structure of a distributive lattice,
as for \emph{Dyck}, \emph{Motzkin} and \emph{Schr\"oder paths}.
Some conditions that guarantee to have a distributive lattice are given in \cite{FP},
where there is also a first attempt to provide a systematic classification of these posets of paths.
Some of these lattices turn out to be isomorphic with the lattices arising from other well known structures
\cite{BarcucciBerniniFerrariPoneti,BF,Ferrari}.
For instance, the lattice generated by all Dyck paths of length $\; n \;$
is isomorphic to the dual of the Young lattice associated with the staircase partition $\; (n,n-1,\ldots,2,1) \,$,
as proved in \cite{S} or in \cite{Sta} where a different combinatorial interpretation is provided.
The language of lattice paths, however, gives a geometric flavor to the subject
which allows to express several properties in a more fascinating way.
Moreover, the study of lattices in terms of paths is a growing area of research,
as witnessed by the investigation of a generalization of Dyck lattices appearing in \cite{Sa}
or by the deep study of the properties of Dyck lattices developed in \cite{STT}.

In the present paper we consider essentially two kinds of questions concerning lattices of paths.
The first topic concerns the representation theory of lattices of \emph{Dyck-like paths},
i.e. paths defined as the ordinary Dyck paths
except for the fact that they use up steps of the form $\; (1,a) \;$ and down steps of the form $\; (1,-b) \,$,
where $\; a \;$ and $\; b \;$ are two assigned positive integers.
Specifically, we obtain a representation theorem for these lattices describing explicitly their spectrum.
We also prove that these lattices can be described as dual of Young lattices associated to a given partition.
Then we review some results scattered in the literature and, in particular, for the Dyck lattices
we observe some elementary properties which will be extensively used in the rest of the paper.

The second topic concerns the \emph{Euler characteristic} of lattices of paths.
The Euler characteristic is a classical invariant measure with play an important role, for instance,
in \emph{combinatorial geometry} \cite{K} and in \emph{geometric probability} \cite{KlainRota}.
The combinatorial interest of the Euler characteristic lies in its deep relation with the M\"{o}bius function \cite{R}
(see also \cite{MunaChar}).
In the present paper, we give a general technique to determine the Euler characteristic of \emph{Dyck-like lattices},
which generalize in a natural way the ordinary Dyck lattices.
In particular, in the case of Dyck and Schr\"{o}der lattices,
we obtain a combinatorial interpretation for the Euler characteristic
in terms of the number of \emph{tunnels} of a path \cite{E1,E2}.
The case of Motzkin lattices cannot be dealt with using the machinery developed for Dyck-like lattices.
So, we provide an \emph{ad hoc} argument which allows to obtain a (not obvious)
combinatorial interpretation of the Euler characteristic also for these lattices.

\section{Basic definitions and properties}

As usual, $\; \NN \;$ denotes the set of all non-negative integers, $\; \ZZ \;$ denotes the set of all integers,
$\; [n] = \{ 1, 2, \ldots, n \} \;$ and $\; [0,n] = \{ 0, 1, \ldots, n \} \,$.
For any real number $\; x \,$, $\; \lfloor x \rfloor \;$ is the greatest integer smaller than $\; x \,$.
In the rest of this section, we will recall some basic definitions concerning lattice paths
and the theory of partially ordered sets \cite{Aigner,Birkhoff,DaveyPriestley}.

\subsection{Lattice paths}
Given a finite subset $\; \Gamma \subseteq \ZZ \,$, a \emph{$\Gamma$-path} of length $\; n \;$
is a function $\; \map{f}{[0,n]}{\NN} \;$ such that $\; f(0) = f(n) = 0 \;$
and $\; f(k+1) - f(k) \in \Gamma \,$, for every $\; k \in [0,n-1] \,$.
Equivalently, a $\Gamma$-path of length $\; n \;$
is a sequence of $\; n \;$ steps $\; (1,k) \,$, with $\; k \in \Gamma \,$,
starting from the origin, ending on the $x$-axis and never going below that axis.

Considered as functions, $\Gamma$-paths of length $\; n \;$ can be ordered coordinatewise
setting $\; f \leq g \;$ whenever $\; f(k) \leq g(k) \;$ for every $\; k \in [0,n] \,$.
In this way, we obtain a poset $\; \CC_n^{\Gamma} \;$ which,
under suitable conditions \cite{FP}, turns out to be a lattice.
In particular, $\; \CC_n^{\Gamma} \;$ is a distributive lattice
when $\; \Gamma = \{ -b, a \} \,$, with $\; a, b \in \NN \,$.
Here $\{-b,a\}$-paths will be called \emph{Dyck-like paths of type} $\; (a,b) \,$,
since they are a natural generalization of ordinary Dyck paths (corresponding to the case $\, a = b = 1 \,$),
and the associated lattice will be denoted by $\; \DD_n^{(a,b)} \;$
(where the subscript $\; n \;$ is related to the length of the paths and its exact meaning will be explained below).
Dyck-like paths have already been considered in \cite{D},
where several results are proved and the case $\; (a,b) = (3,2) \;$ is examined in great detail,
and also in \cite{BaFl}, where they are used as a source of examples.

We will use $\; U \,$, $\; H \;$ and $\; D \;$ to denote up steps, horizontal steps and down steps.
Given a path $\; \gamma \;$ having precisely one type of up step and one type of down step,
a \emph{peak} of $\; \gamma \;$ is just a sequence formed by an up step and a down step.
An \emph{elevated Dyck path} is a Dyck path touching the $x$-axis just at its starting and ending points.

In the sequel, a path will be denoted by a Latin letter when considered as an element of a lattice
and by a Greek letter in all other cases.

\subsection{Partially ordered sets}
In a poset $\; P \,$, an element $\; x \;$ is \emph{covered} by an element $\; y \;$
when $\; x \leq u \leq y \;$ implies $\; x = u \;$ or $\; u = y \,$.
A poset $\; P \;$ is \emph{ranked} when it admits a \emph{rank function},
that is a function $\; \map{r}{P}{\NN} \;$ such that
$\; r(y) = r(x) + 1 \;$ whenever $\; x \;$ is covered by $\; y \,$.
For finite posets the function $\; r \;$ is usually chosen so that the minimal elements have rank $\; 0 \,$.
The \emph{height} of $\; P \;$ is its maximum rank.

A \emph{join-semilattice} (\emph{meet-semilattice}) is a poset in which
there exists the supremum (the infimum) of any two elements.
A lattice is a poset in which there exists the supremum and the infimum of any two elements.

A poset $\; P \;$ has a \emph{minimum} $\; \widehat{0} \;$ (\emph{maximum} $\; \widehat{1} \,$)
when it has only one minimal (maximal) element.
In a poset $\; P \;$ with minimum (maximum),
an \emph{atom} (\emph{coatom}) is an element covering the minimum (covered by the maximum).
In a (finite) lattice $\; L \,$, the \emph{socle} is the join of all atoms
and the \emph{radical} is the meet of all coatoms.

An \emph{order-ideal} of a poset $\; P \;$ is a subset $\; I \;$ such that
$\; x \in I \;$ and $\; u \leq x \;$ imply $\; u \in I \,$.
The \emph{principal ideal} $\; \downarrow x \;$ generated by an element $\; x \in P \;$
is the set of all elements $\; u \in P \;$ such that $\; u \leq x \,$.
Similarly, the \emph{principal filter} $\; \uparrow x \;$ generated by an element $\; x \in P \;$
is the set of all elements $\; u \in P \;$ such that $\; u \geq x \,$.

A \emph{join-irreducible element} of a distributive lattice $\; \DD \;$
is any element $\; x \ne \widehat{0} \;$ with the property that
if $\; x = u \vee v \;$ then $\; x = u \;$ or $\; x = v \,$.
The set $\; \JJ(P) \;$ of all order-ideals of $\; P \,$, ordered by inclusion, is a distributive lattice.
Conversely, by Birkhoff's representation theorem,
every finite distributive lattice $\; \DD \;$ is isomorphic to the lattice $\; \JJ(P) \;$
where $\; P = \Spec(\DD) \,$, the \emph{spectrum} of $\; \DD \,$,
is defined as the poset of all join-irreducibles of $\; \DD \,$.

A \emph{valuation} on a distributive lattice $\; \DD \;$ with values in $\; \RR \;$
is a function $\; \map{\nu}{\DD}{\RR}$ such that $\; \nu(\widehat{0}) = 0 \;$ and
$\; \nu( x \vee y ) + \nu( x \wedge y ) = \nu(x) + \nu(y) \;$ for every $\; x, y \in \DD \,$.
A valuation on a finite distributive lattice $\; \DD \;$
is uniquely determined by the values it takes on the set of join-irreducibles of $\; \DD \,$,
and these values can be arbitrarily assigned \cite{R}.
Every valuation $\; \nu \;$ satisfies the following generalized form of the \emph{principle of inclusion-exclusion}:
\begin{equation}\label{dasilva}
 \nu( x_1 \vee \cdots \vee x_n ) =
 \sum_{S \subseteq [n] \atop S \neq \emptyset } (-1)^{|S|-1} \; \nu\!\left( \bigwedge_{i\in S}x_i \right)\, .
\end{equation}

The \emph{Euler characteristic} of $\; \DD \;$ is defined as the unique valuation $\; \chi \;$
such that $\; \chi( \widehat{0} ) = 0 $ and $\; \chi(x) = 1 \;$ for every join-irreducible $\; x \;$ of $\; \DD \,$.
In particular, $\; \chi(D) = \chi(\widehat{1}) \,$.

A map $\; \map{f}{P}{Q} \;$ between two posets $\; P \;$ and $\; Q \;$
is \emph{order-preserving} when $\; x \leq y \;$ implies $\; f(x) \leq f(y) \;$ for every $\; x, y \in P \,$.
In particular, it is a \emph{poset isomorphism} when it is an order-preserving bijection.
In this paper, all isomorphisms between posets are always understood to be poset isomorphisms.

A \emph{partition} of a non-negative integer $\; n \;$ is a sequence $\; \lambda = (\ll_1,\ll_2,\ldots,\ll_k) \;$
of positive integers such that $\; \ll_1 \geq \ll_2 \geq \cdots \geq \ll_k > 0 \;$
and $\; \ll_1 + \ll_2 + \cdots + \ll_k = n \,$.
The $\; \ll_i$'s are called the \emph{parts} of $\; \ll \,$,
and the sum of all parts of $\; \ll \;$ will be denoted by $\; |\lambda | \,$.
The (\emph{Ferrers}) \emph{diagram} of $\; \ll \;$ is a left-justified array of squares (or dots)
with exactly $\; \ll_i \;$ squares in the $\; i$-th row.
Partitions can be ordered by magnitude of parts \cite{Aigner}:
if $\; \alpha = (a_1,\ldots,a_h) \;$ and $\; \beta = (b_1,\ldots,b_k) \,$, then $\; \alpha \leq \beta \;$
whenever $\; h \leq k \;$ and $\; a_i \leq b_i \;$ for every $\; i = 1, 2, \ldots, h \,$.
If $\; \alpha \leq \beta \;$ the diagram of $\; \alpha \;$ is contained in the diagram of $\; \beta \,$.
The resulting poset is an infinite distributive lattice, called \emph{Young lattice}.
In particular, the \emph{Young lattice} $\; \YY_\lambda \;$ generated by a partition $\; \lambda \;$
is the set of all integer partitions $\; \alpha \;$ such that $\; \alpha \leq \lambda \,$,
i.e. $\; \YY_\lambda =\; \downarrow \ll \;$ in $\; \YY \,$.

\section{Dyck-like paths}

\subsection{Representation of Dyck-like paths}\label{sec2}

To study Dyck-like paths of type $\; (a,b) \;$ we can always suppose, without loss of generality, that $\; a \geq b \,$.
We begin our study noticing that the length of a Dyck-like path of type $\; (a,b) \;$
strictly depends on $\; a \;$ and $\; b \,$, as stated in the following proposition essentially due to Duchon \cite{D}.
\begin{proposition}\label{lengthDyckLikePaths}
 Any Dyck-like path of type $\; (a,b) \;$ starting from the origin
 ends at the point $\; (n\cdot \ell(a,b), 0 ) \,$,
 where $\; \ell(a,b) = \frac{a+b}{\gcd(a,b)} \;$ and $\; n \in \NN \,$.
\end{proposition}
\begin{proof}
 First consider a Dyck-like path of type $\; (a,b) \;$ of minimum length
 (i.e. consisting of the minimum positive number of steps), belonging to $\; \DD_1^{(a,b)} \,$.
 It is made of $\; h \;$ steps $\; (1,a) \;$ and $\; k \;$ steps $\; (1,-b) \,$,
 where $\; h \;$ and $\; k \;$ are the minimum positive integers such that $\; h a - k b = 0 \,$.
 Since $\; h = \frac{b}{\gcd(a,b)} \;$ and $\; k = \frac{a}{\gcd(a,b)} \,$,
 it follows that it has length $\; \ell(a,b) = \frac{a+b}{\gcd(a,b)} \,$.
 The length of any other Dyck-like path of type $\; (a,b) \;$ is a multiple of $\; \ell(a,b) \,$.
\end{proof}

\begin{figure}[h]
\begin{center}
 \setlength{\unitlength}{3mm}
 \begin{picture}(26,8)
 \put(0,0){
 \begin{picture}(5,4)(0,-1)
  \multiput(0,0)(1,0){5}{\dashbox{0.05}(1,4){}}
  \multiput(0,0)(0,1){4}{\dashbox{0.05}(5,1){}}
  \put(0,0){\circle*{0.3}}
  \put(1,3){\circle*{0.3}}
  \put(2,1){\circle*{0.3}}
  \put(3,4){\circle*{0.3}}
  \put(4,2){\circle*{0.3}}
  \put(5,0){\circle*{0.3}}
  \put(0,0){\line(1,3){1}}
  \put(1,3){\line(1,-2){1}}
  \put(2,1){\line(1,3){1}}
  \put(3,4){\line(1,-2){2}}
  \put(2.5,-1){\makebox(0,0){(i)}}
 \end{picture}}
 \put(8,0){
 \begin{picture}(7,7)(0,-1)
  \multiput(0,0)(1,0){7}{\dashbox{0.05}(1,6){}}
  \multiput(0,0)(0,1){6}{\dashbox{0.05}(7,1){}}
  \put(0,0){\circle*{0.3}}
  \put(1,5){\circle*{0.3}}
  \put(2,3){\circle*{0.3}}
  \put(3,1){\circle*{0.3}}
  \put(4,6){\circle*{0.3}}
  \put(5,4){\circle*{0.3}}
  \put(6,2){\circle*{0.3}}
  \put(7,0){\circle*{0.3}}
  \put(0,0){\line(1,5){1}}
  \put(1,5){\line(1,-2){2}}
  \put(3,1){\line(1,5){1}}
  \put(4,6){\line(1,-2){3}}
  \put(3.5,-1){\makebox(0,0){(ii)}}
 \end{picture}}
 \put(18,0){
 \begin{picture}(8,7)(0,-1)
  \multiput(0,0)(1,0){8}{\dashbox{0.05}(1,7){}}
  \multiput(0,0)(0,1){7}{\dashbox{0.05}(8,1){}}
  \put(0,0){\circle*{0.3}}
  \put(1,5){\circle*{0.3}}
  \put(2,2){\circle*{0.3}}
  \put(3,7){\circle*{0.3}}
  \put(4,4){\circle*{0.3}}
  \put(5,1){\circle*{0.3}}
  \put(6,6){\circle*{0.3}}
  \put(7,3){\circle*{0.3}}
  \put(8,0){\circle*{0.3}}
  \put(0,0){\line(1,5){1}}
  \put(1,5){\line(1,-3){1}}
  \put(2,2){\line(1,5){1}}
  \put(3,7){\line(1,-3){2}}
  \put(5,1){\line(1,5){1}}
  \put(6,6){\line(1,-3){2}}
  \put(4,-1){\makebox(0,0){(iii)}}
 \end{picture}}
 \end{picture}
\end{center}
 \caption{The minimum element in $\; \mathcal{\DD}_1^{(a,b)}$ for
 (i) $\; a = 3 \,$, $\; b = 2 \,$, (ii) $\; a = 5 \,$, $\; b = 2 \;$ and (iii) $\; a = 5 \,$, $\; b = 3 \,$.}
 \label{figuraSP}
\end{figure}
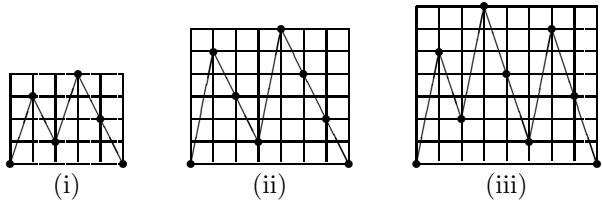

The maximum of $\; \DD_1^{(a,b)} \;$ is the path consisting of $\; \frac{b}{\gcd(a,b)} \;$ steps $\; (1,a) \;$
followed by $\; \frac{a}{\gcd(a,b)} \;$ steps $\; (1,-b) \,$,
whereas the minimum is obtained by starting with a step $\; (1,a) \;$
followed by as many steps $\; (1,-b) \;$ as possible (i.e. without going below the $x$-axis),
and then repeating this procedure until we reach for the first time the $x$-axis
(see Figure \ref{figuraSP} for some examples).

Our next goal is to find some suitable representation results for the lattices $\; \DD_n^{(a,b)} \;$
of all Dyck-like paths of type $\; (a,b) \;$ having length $\; n \cdot \ell(a,b) \,$.
We begin by giving an alternative description of $\; \DD_n^{(a,b)} \,$.
For each path in $\; \DD_n^{(a,b)} \,$, consider the path obtained by replacing
each step $\; (1,a) \;$ with a step $\; (1,1) \;$ and each step $\; (1,-b) \;$ with a step $\; (1,-1) \,$.
If $\; a \neq b \,$, the paths of the resulting set $\; \overline{\DD}_n^{(a,b)}$ terminates below the $x$-axis.
\begin{lemma}
 The paths in $\; \overline{\DD}_n^{(a,b)} \;$ start from the origin,
 end at the point $\; ( n \cdot \frac{a+b}{\gcd(a,b)}, n \cdot \frac{b-a}{\gcd(a,b)} ) \;$
 and never go below the line $\; y = \frac{b-a}{a+b} x \,$.
\end{lemma}
\begin{proof}
 If $\; (x_0,y_0) \;$ is the endpoint of the paths in $\; \overline{\DD}_n^{(a,b)} \;$
 then $\; x_0 \;$ is equal to the length of the paths in $\; \DD_n^{(a,b)} \;$
 obtained in Proposition \ref{lengthDyckLikePaths},
 and $\; y_0 \;$ can be obtained by subtracting from the number $\; h \;$ of steps $\; (1,1) \;$
 the number $\; k \;$ of steps $\; (1,-1) \,$.
 Since $\; h \;$ and $\; k \;$ are defined as in Proposition \ref{lengthDyckLikePaths},
 it follows that $\; y_0 = n \cdot \frac{b-a}{\gcd(a,b)} \,$.
 The last part of the thesis is obtained by determining the line passing through the origin and the point $\; (x_0,y_0) \,$.
\end{proof}

\begin{figure}[h]
\begin{center}
 \setlength{\unitlength}{3mm}
 \begin{picture}(47,16)
 \put(0,0){
 \setlength{\unitlength}{2.5mm}
 \begin{picture}(15,18)
  \put(0,0){\vector(1,0){17}}
  \put(0,0){\vector(0,1){19}}
  \multiput(0,0)(1,3){7}{\circle*{0.3}}
  \multiput(2,1)(1,3){6}{\circle*{0.3}}
  \multiput(4,2)(1,3){5}{\circle*{0.3}}
  \multiput(5,0)(1,3){5}{\circle*{0.3}}
  \multiput(7,1)(1,3){4}{\circle*{0.3}}
  \multiput(9,2)(1,3){3}{\circle*{0.3}}
  \multiput(10,0)(1,3){3}{\circle*{0.3}}
  \multiput(12,1)(1,3){2}{\circle*{0.3}}
  \put(14,2){\circle*{0.3}}
  \put(15,0){\circle*{0.3}}
  \put(0,0){\line(1,3){6}}
  \put(2,1){\line(1,3){5}}
  \put(4,2){\line(1,3){4}}
  \put(5,0){\line(1,3){4}}
  \put(7,1){\line(1,3){3}}
  \put(9,2){\line(1,3){2}}
  \put(10,0){\line(1,3){2}}
  \put(12,1){\line(1,3){1}}
  \put(2,1){\line(-1,2){1}}
  \put(5,0){\line(-1,2){3}}
  \put(7,1){\line(-1,2){4}}
  \put(10,0){\line(-1,2){6}}
  \put(12,1){\line(-1,2){7}}
  \put(15,0){\line(-1,2){9}}
 \end{picture}}
 \qbezier(12,7)(13,8.7)(14,9)
 \put(14,9){\vector(4,1){0.3}}
 \put(14.5,3){
 \setlength{\unitlength}{3.6mm}
 \begin{picture}(15,10)(0,-3.5)
  \put(0,0){\vector(1,0){15}}
  \put(0,-3.5){\vector(0,1){10}}
  \multiput(0,0)(1,1){7}{\circle*{0.2}}
  \multiput(2,0)(1,1){6}{\circle*{0.2}}
  \multiput(4,0)(1,1){5}{\circle*{0.2}}
  \multiput(5,-1)(1,1){5}{\circle*{0.2}}
  \multiput(7,-1)(1,1){4}{\circle*{0.2}}
  \multiput(9,-1)(1,1){3}{\circle*{0.2}}
  \multiput(10,-2)(1,1){3}{\circle*{0.2}}
  \multiput(12,-2)(1,1){2}{\circle*{0.2}}
  \put(14,-2){\circle*{0.2}}
  \put(15,-3){\circle*{0.2}}
  \put(0,0){\line(1,1){6}}
  \put(2,0){\line(1,1){5}}
  \put(4,0){\line(1,1){4}}
  \put(5,-1){\line(1,1){4}}
  \put(7,-1){\line(1,1){3}}
  \put(9,-1){\line(1,1){2}}
  \put(10,-2){\line(1,1){2}}
  \put(12,-2){\line(1,1){1}}
  \put(2,0){\line(-1,1){1}}
  \put(5,-1){\line(-1,1){3}}
  \put(7,-1){\line(-1,1){4}}
  \put(10,-2){\line(-1,1){6}}
  \put(12,-2){\line(-1,1){7}}
  \put(15,-3){\line(-1,1){9}}
  \put(0,0){\line(5,-1){15}}
 \end{picture}}
 \qbezier(33,9)(34,9.9)(35,10)
 \put(35,10){\vector(4,1){0.3}}
 \put(36,8){
 \setlength{\unitlength}{3.4mm}
 \begin{picture}(9,6)
  \put(0,0){\dashbox{0.2}(9,6){}}
  \multiput(0,0)(0,1){7}{\circle*{0.2}}
  \multiput(1,1)(0,1){6}{\circle*{0.2}}
  \multiput(2,2)(0,1){5}{\circle*{0.2}}
  \multiput(3,2)(0,1){5}{\circle*{0.2}}
  \multiput(4,3)(0,1){4}{\circle*{0.2}}
  \multiput(5,4)(0,1){3}{\circle*{0.2}}
  \multiput(6,4)(0,1){3}{\circle*{0.2}}
  \multiput(7,5)(0,1){2}{\circle*{0.2}}
  \put(8,6){\circle*{0.2}}
  \put(9,6){\circle*{0.2}}
  \put(0,0){\line(0,1){6}}
  \put(1,1){\line(0,1){5}}
  \put(2,2){\line(0,1){4}}
  \put(3,2){\line(0,1){4}}
  \put(4,3){\line(0,1){3}}
  \put(5,4){\line(0,1){2}}
  \put(6,4){\line(0,1){2}}
  \put(7,5){\line(0,1){1}}
  \put(0,1){\line(1,0){1}}
  \put(0,2){\line(1,0){3}}
  \put(0,3){\line(1,0){4}}
  \put(0,4){\line(1,0){6}}
  \put(0,5){\line(1,0){7}}
  \put(0,6){\line(1,0){9}}
 \end{picture}}
 \put(41,6.7){\makebox(0,0){$\downarrow$}}
 \put(36,0){
 \setlength{\unitlength}{3.4mm}
 \begin{picture}(7,5)
  \put(0,0){\line(0,1){5}}
  \put(1,0){\line(0,1){5}}
  \put(2,1){\line(0,1){4}}
  \put(3,1){\line(0,1){4}}
  \put(4,2){\line(0,1){3}}
  \put(5,3){\line(0,1){2}}
  \put(6,3){\line(0,1){2}}
  \put(7,4){\line(0,1){1}}
  \put(0,0){\line(1,0){1}}
  \put(0,1){\line(1,0){3}}
  \put(0,2){\line(1,0){4}}
  \put(0,3){\line(1,0){6}}
  \put(0,4){\line(1,0){7}}
  \put(0,5){\line(1,0){7}}
  \put(0.2,1.2){\dashbox{0.2}(2.6,1.6){}}
  \put(0.2,3.2){\dashbox{0.2}(2.6,1.6){}}
  \put(3.2,3.2){\dashbox{0.2}(2.6,1.6){}}
 \end{picture}}
 \end{picture}
\end{center}
 \caption{A path in $\; \DD_n^{(a,b)} \;$ and the associated partition for $\; a = 3 \,$, $\; b = 2 \;$ and $\; n = 3 \,$.}
 \label{PartitionsFig}
\end{figure}
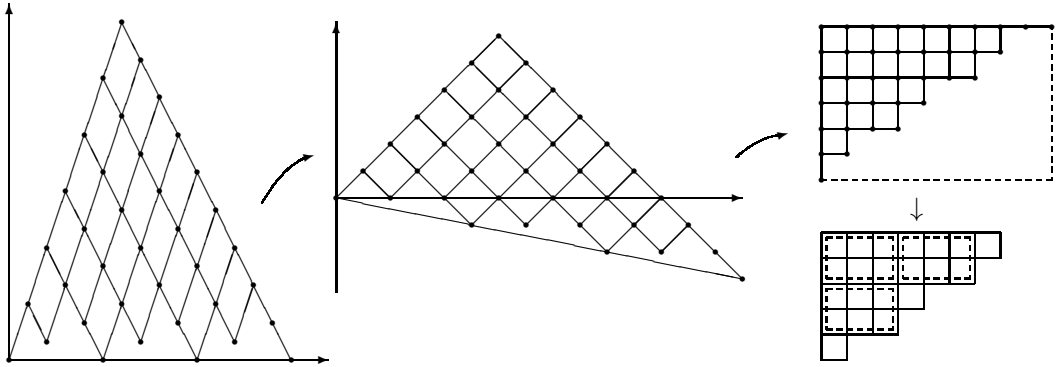

A careful inspection of some examples (as the one in Figure \ref{PartitionsFig})
shows that the paths in $\; \overline{\DD}_{n}^{(a,b)} \;$ (or, equivalently, in $\; \DD_n^{(a,b)} \;$)
are in bijection with the integer partitions whose Ferrers diagram is included
in the Ferrers diagram obtained by taking simultaneously
the minimum and the maximum paths in $\; \overline{\DD}_n^{(a,b)} \,$,
then rotating by $\; 45^{\circ} \;$ anticlockwise
and finally considering the squares obtained by drawing that part of the lattice grid included between the two paths.
We will denote such a partition $\; \lambda_n^{(a,b)} \,$,
and call it the \emph{partition associated with} $\; \overline{\DD}_n^{(a,b)} \;$ (or $\; \DD_n^{(a,b)} \,$).

All this implies
\begin{proposition}\label{young}
 The distributive lattice $\; \DD_{n}^{(a,b)} \;$ is isomorphic
 to the dual of the Young lattice $\; \YY_{\lambda_{n}^{(a,b)}} \,$.
\end{proposition}

Before proceeding further,
it will be useful to characterize the join-irreducible elements of $\; \DD_n^{(a,b)} \,$.
It is easy to see that a path in $\; \DD_n^{(a,b)} \;$ is join-irreducible
if and only if there is precisely one peak $\; (1,a)(1,-b) \;$
which can be replaced by the valley $\; (1,-b)(1,a) \;$ without letting the path cross the $x$-axis.
If we interpret a path as a word on the $2$-letter alphabet $\; \{ U, D \} \,$,
where $\; U = (1,a) \;$ and $\; D = (1,-b) \,$,
then the above condition can be easily translated as follows.
A path $\; w \in \DD_n^{(a,b)} \;$ is join-irreducible if and only if
it is obtained from the minimum element $\; \widehat{0} \;$ of the lattice
by replacing a subword $\; w' \;$ of length $\; k \geq 2 \;$
starting with $\; D \;$ and ending with $\; U \;$ with the word $\; U^r D^s \,$,
for suitable $\; r \;$ and $\; s \;$ such that $\; r + s = k \,$.
Each maximal word of type $\; U^r D^s \;$ appearing as a subword of a path $\; w \;$
will be called a \emph{pyramid} of $\; w \;$ when, replacing it with $\; U^{r-1} D U D^{s-1} \,$,
the resulting path is still in $\; \DD_n^{(a,b)} \,$.
In conclusion, we can say that a path of $\; \DD_n^{(a,b)} \;$ is join-irreducible if and only if it has a unique pyramid.
Of course, all what we have said concerning the join-irreducibles of $\; \DD_n^{(a,b)} \;$
can be easily transferred to $\; \overline{\DD}_n^{(a,b)} \,$.

Now, we can state our main result concerning the representation of Dyck-like paths.
\begin{theorem}\label{rep}
 The distributive lattice $\; \DD_n^{(a,b)} \;$ is isomorphic to $\; \JJ(P_{n}^{(a,b)})$,
 where $\; P_{n}^{(a,b)} \;$ is the set of points $\; (x,y)\in [0,n\cdot \ell (a,b)]\times \ZZ \;$
 such that $\; y \leq x \,$, $\; y \leq -x+n\frac{2b}{\gcd (a,b)} \,$, $\; y\geq \frac{b-a}{a+b} \,$,
 ordered coordinatewise.
\end{theorem}
\begin{proof}
 Proposition \ref{young} suggests that the set of join-irreducibles of $\; \overline{\DD}_n^{(a,b)} \;$
 is in bijection with the cells of the Ferrers diagram of $\; \lambda_n^{(a,b)} \;$ (see Figure \ref{PartitionsFig}).
 Indeed, we can associate with each join-irreducible $\; x \;$ the cell containing
 the vertex of the unique pyramid of $\; x \;$
 which is included between $\; x \;$ and the line $\; y = \frac{b-a}{a+b} x \,$.
 Therefore, identifying the cells of the Ferrers diagram with their topmost vertex,
 we obtain a bijection (actually, an isomorphism)
 between the set of join-irreducibles of $\; \overline{\DD}_n^{(a,b)} \;$
 (and so of $\; \DD_n^{(a,b)} \,$) and the set $\; P_n^{(a,b)} \,$.
 Now, the thesis follows from Birkhoff's representation theorem for finite distributive lattices.
\end{proof}

Next proposition describes the partition $\; \lambda_n^{(a,b)} \;$
in terms of the positive integers $\; a \,$, $\; b \;$ and $\; n \,$,
when $\; \gcd(a,b) = 1 \;$ (the general case follows immediately).
\begin{proposition}\label{partition}
 Let $\; \gcd(a,b) = 1 \,$. For $\; n = 1 \,$,
 \begin{equation}\label{partition1}
  \lambda_1^{(a,b)} =
  \left(
   \left\lfloor (b-1)\cdot \frac{a}{b}\right\rfloor,
   \left\lfloor (b-2)\cdot \frac{a}{b}\right\rfloor,
   \ldots,
   \left\lfloor 2\cdot \frac{a}{b}\right\rfloor,
   \left\lfloor \frac{a}{b}\right\rfloor
  \right)
 \end{equation}
 is a partition with $\; b - 1 \;$ parts.
 More generally, for any $\; n \in \NN \,$,
 $$
  \lambda_{n}^{(a,b)} =
  ( \ll_{n,b}, \ldots, \ll_{n,1}, \ll_{n-1,b}, \ldots, \ll_{n-1,1}, \ldots,
   \ll_{2,b}, \ldots, \ll_{2,1}, \ll_{1,b}, \ldots, \ll_{1,2} )
 $$
 is a partition with $\; n b - 1 \;$ parts, where
 $\; \ll_{h,k} = (h-1) a + \left\lfloor (k-1) \cdot a/b \right\rfloor \,$.
 In particular, $$ |\lambda_n^{(a,b)}| = a b \frac{n(n-1)}{2} + n |\lambda_1^{(a,b)}| \, . $$
\end{proposition}
\begin{proof}
 For $\; n = 1 \,$, let $\; \lambda_1^{(a,b)} = (k_{b-1},\ldots,k_2,k_1) \;$
 with $\; k_{b-1} \geq \cdots \geq k_2 \geq k_1 \,$.
 From the form of the paths in $\; \overline{\DD}_n^{(a,b)} \;$ and the definition of $\; \lambda_1^{(a,b)} \,$,
 it follows that $\; k_i \;$ is the sum of the cardinalities
 of the first $\; i \;$ sequences of consecutive down steps of the minimum path of $\; \DD_n^{(a,b)} \,$.
 From the definition of this minimum, it follows that $\; k_i \;$ is defined by the inequalities
 $\; i a - k_i b \geq 0 \;$ and $\; i a - (k_i +1) b < 0 \,$,
 or equivalently by $\; k_i \leq i\; a/b < k_i + 1 \,$,
 and hence $\; k_i = \left\lfloor i\cdot a/b \right\rfloor \,$.

 For an arbitrary $\; n \,$,
 just observe that the Ferrers diagram of $\; \lambda_n^{(a,b)} \;$
 is a staircase-like diagram made of $\; (a \times b)$-rectangles,
 where the topmost row consists of $\; n - 1 \;$ rectangles and,
 at the end of each horizontal strip of rectangles,
 the Ferrers diagram of $\; \lambda_1^{(a,b)} \;$ is appended (see Figure \ref{PartitionsFig}).
\end{proof}

The rank of a path $\; x \in \DD_n^{(a,b)} \;$ can be easily expressed in terms of its \emph{area} $\; \AA(x) \,$,
i.e. the area of the region determined by the path and the $x$-axis.
\begin{proposition}
 The rank of any element $\; x \in \DD_n^{(a,b)} \;$ is
 \begin{equation}\label{rankDn}
  r(x) = \frac{\AA(x)-\AA(\widehat{0})}{a+b} \, .
 \end{equation}
\end{proposition}
\begin{proof}
 We will prove that the function defined in (\ref{rankDn}) satisfies the properties of a rank function.
 First we have $\; r(\widehat{0}) = \frac{\AA(\widehat{0})-\AA(\widehat{0})}{a+b} = 0 \,$.
 Suppose now that $\; x \;$ is covered by $\; y \,$.
 Then the path $\; y \;$ is obtained from $\; x \;$ by replacing a valley with a peak.
 This implies that the area $\; A(y) \;$ is obtained from $\; \AA(x) \;$
 by adding the area of a parallelogram of area $\; a + b \,$.
 Hence $\; r(y) = \frac{\AA(x)+a+b-\AA(\widehat{0})}{a+b} = r(x) + 1 \,$.
\end{proof}

\subsection{Dyck lattices}

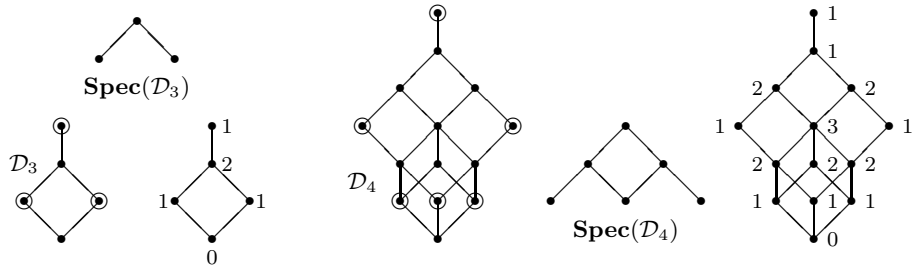
\begin{figure}[h]
\begin{center}
 \setlength{\unitlength}{5mm}
 \begin{picture}(23,6)
   \put(0,0){
   \begin{picture}(2,4)
    \put(1,0){\circle*{0.2}}
    \put(0,1){\circle*{0.2}} \put(0,1){\circle{0.4}}
    \put(2,1){\circle*{0.2}} \put(2,1){\circle{0.4}}
    \put(1,2){\circle*{0.2}}
    \put(1,3){\circle*{0.2}} \put(1,3){\circle{0.4}}
    \put(1,0){\line(1,1){1}}
    \put(1,0){\line(-1,1){1}}
    \put(1,2){\line(-1,-1){1}}
    \put(1,2){\line(1,-1){1}}
    \put(1,2){\line(0,1){1}}
    \put(0,2){\makebox(0,0){\small$\DD_3$}}
   \end{picture}}
  \put(4,0){
   \begin{picture}(2,4)
    \put(1,0){\circle*{0.2}} \put(1,-0.5){\makebox(0,0){\footnotesize$0$}}
    \put(0,1){\circle*{0.2}} \put(-0.3,1){\makebox(0,0){\footnotesize$1$}}
    \put(2,1){\circle*{0.2}} \put(2.3,1){\makebox(0,0){\footnotesize$1$}}
    \put(1,2){\circle*{0.2}} \put(1.4,2){\makebox(0,0){\footnotesize$2$}}
    \put(1,3){\circle*{0.2}} \put(1.4,3){\makebox(0,0){\footnotesize$1$}}
    \put(1,0){\line(1,1){1}}
    \put(1,0){\line(-1,1){1}}
    \put(1,2){\line(-1,-1){1}}
    \put(1,2){\line(1,-1){1}}
    \put(1,2){\line(0,1){1}}
   \end{picture}}
  \put(2,4.8){
   \begin{picture}(2,4)
    \put(0,0){\circle*{0.2}}
    \put(2,0){\circle*{0.2}}
    \put(1,1){\circle*{0.2}}
    \put(0,0){\line(1,1){1}}
    \put(2,0){\line(-1,1){1}}
    \put(1,-0.8){\makebox(0,0){\small$\Spec(\DD_3)$}}
   \end{picture}}
  \put(9,0){
   \begin{picture}(4,6)
    \put(2,0){\circle*{0.2}}
    \put(1,1){\circle*{0.2}} \put(1,1){\circle{0.4}}
    \put(2,1){\circle*{0.2}} \put(2,1){\circle{0.4}}
    \put(3,1){\circle*{0.2}} \put(3,1){\circle{0.4}}
    \put(1,2){\circle*{0.2}}
    \put(2,2){\circle*{0.2}}
    \put(3,2){\circle*{0.2}}
    \put(0,3){\circle*{0.2}} \put(0,3){\circle{0.4}}
    \put(2,3){\circle*{0.2}}
    \put(4,3){\circle*{0.2}} \put(4,3){\circle{0.4}}
    \put(1,4){\circle*{0.2}}
    \put(3,4){\circle*{0.2}}
    \put(2,5){\circle*{0.2}}
    \put(2,6){\circle*{0.2}} \put(2,6){\circle{0.4}}
    \put(2,0){\line(-1,1){1}}
    \put(2,0){\line(0,1){1}}
    \put(2,0){\line(1,1){1}}
    \put(2,1){\line(-1,1){2}}
    \put(2,1){\line(1,1){2}}
    \put(1,2){\line(1,1){2}}
    \put(3,2){\line(-1,1){2}}
    \put(0,3){\line(1,1){2}}
    \put(4,3){\line(-1,1){2}}
    \put(1,1){\line(0,1){1}}
    \put(3,1){\line(0,1){1}}
    \put(2,2){\line(0,1){1}}
    \put(2,5){\line(0,1){1}}
    \put(2,2){\line(-1,-1){1}}
    \put(2,2){\line(1,-1){1}}
    \put(0,1.5){\makebox(0,0){\small$\DD_4$}}
   \end{picture}}
  \put(14,1){
   \begin{picture}(4,2)
    \put(0,0){\circle*{0.2}}
    \put(2,0){\circle*{0.2}}
    \put(4,0){\circle*{0.2}}
    \put(1,1){\circle*{0.2}}
    \put(3,1){\circle*{0.2}}
    \put(2,2){\circle*{0.2}}
    \put(0,0){\line(1,1){2}}
    \put(4,0){\line(-1,1){2}}
    \put(1,1){\line(1,-1){1}}
    \put(3,1){\line(-1,-1){1}}
    \put(2,-0.8){\makebox(0,0){\small$\Spec(\DD_4)$}}
   \end{picture}}
  \put(19,0){
   \begin{picture}(4,6)
    \put(2,0){\circle*{0.2}} \put(2.5,0){\makebox(0,0){\footnotesize$0$}}
    \put(1,1){\circle*{0.2}} \put(0.5,1){\makebox(0,0){\footnotesize$1$}}
    \put(2,1){\circle*{0.2}} \put(2.5,1){\makebox(0,0){\footnotesize$1$}}
    \put(3,1){\circle*{0.2}} \put(3.5,1){\makebox(0,0){\footnotesize$1$}}
    \put(1,2){\circle*{0.2}} \put(0.5,2){\makebox(0,0){\footnotesize$2$}}
    \put(2,2){\circle*{0.2}} \put(2.5,2){\makebox(0,0){\footnotesize$2$}}
    \put(3,2){\circle*{0.2}} \put(3.5,2){\makebox(0,0){\footnotesize$2$}}
    \put(0,3){\circle*{0.2}} \put(-0.5,3){\makebox(0,0){\footnotesize$1$}}
    \put(2,3){\circle*{0.2}} \put(2.5,3){\makebox(0,0){\footnotesize$3$}}
    \put(4,3){\circle*{0.2}} \put(4.5,3){\makebox(0,0){\footnotesize$1$}}
    \put(1,4){\circle*{0.2}} \put(0.5,4){\makebox(0,0){\footnotesize$2$}}
    \put(3,4){\circle*{0.2}} \put(3.5,4){\makebox(0,0){\footnotesize$2$}}
    \put(2,5){\circle*{0.2}} \put(2.5,5){\makebox(0,0){\footnotesize$1$}}
    \put(2,6){\circle*{0.2}} \put(2.5,6){\makebox(0,0){\footnotesize$1$}}
    \put(2,0){\line(-1,1){1}}
    \put(2,0){\line(0,1){1}}
    \put(2,0){\line(1,1){1}}
    \put(2,1){\line(-1,1){2}}
    \put(2,1){\line(1,1){2}}
    \put(1,2){\line(1,1){2}}
    \put(3,2){\line(-1,1){2}}
    \put(0,3){\line(1,1){2}}
    \put(4,3){\line(-1,1){2}}
    \put(1,1){\line(0,1){1}}
    \put(3,1){\line(0,1){1}}
    \put(2,2){\line(0,1){1}}
    \put(2,5){\line(0,1){1}}
    \put(2,2){\line(-1,-1){1}}
    \put(2,2){\line(1,-1){1}}
   \end{picture}}
 \end{picture}
\end{center}
\caption{The Dyck lattices $\; \DD_3 \;$ and $\; \DD_4 \,$, their spectra and the distribution of the characteristic.}
\label{DyckLatticeD6}
\end{figure}

In this section we will consider the special case of ordinary Dyck lattices $\; \DD_n \,$,
corresponding to the case $\; (a,b) = (1,1) \,$.
$\; \DD_n \;$ is a distributive lattice of height $\; { n \choose 2 } \,$,
with minimum $\; (UD)^n \;$ and maximum $\; U^n D^n \,$,
with $\; n-1 \;$ atoms of the form $\; (UD)^kUUDD(UD)^{2n-4k} \;$ and just $\; 1 \;$ coatom $\; U^{n-1}DUD^{n-1} \,$.
It is easy to see that the socle of $\; \DD_n \;$ is the path $\; s = U(UD)^{n-1}D \;$
and that the principal ideal $\; \downarrow s \;$ is isomorphic to a Boolean algebra $\; B_{n-1} \;$
while the principal filter $\; \uparrow s \;$ is isomorphic to a Dyck lattice $\; \DD_{n-1} \,$, whenever $\; n \geq 1 \;$
(see Figure \ref{DyckLatticeD6}).
The rank function can be expressed in terms of the area, namely $\; r(x) = (\AA(x)-n)/2 \,$.

A \emph{pyramid} in a Dyck path is a maximal sequence of consecutive steps of the form $\; U^h D^h \,$,
for some $\; h \geq 1 \,$,
which can be replaced with $\; U^{h-1} D U D^{h-1} \;$ still remaining inside the class of Dyck paths.
The positive integer $\; h \;$ is called the \emph{dimension} of the pyramid,
whereas the height of the vertex is called the \emph{height} of the pyramid.
The present definition of pyramid for Dyck paths
is a special case of the definition given in section \ref{sec2} for Dyck-like paths.
\begin{proposition}\label{SpecDn}
 The join-irreducibles of the lattice $\; \DD_n \;$ are the paths with exactly one pyramid,
 that is the paths of the form
 \begin{center}
 \setlength{\unitlength}{4mm}
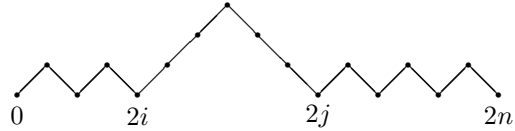
 \begin{picture}(16,4)(0,-1)
  \put(0,0){\circle*{0.2}}
  \put(1,1){\circle*{0.2}}
  \put(2,0){\circle*{0.2}}
  \put(3,1){\circle*{0.2}}
  \put(4,0){\circle*{0.2}}
  \put(5,1){\circle*{0.2}}
  \put(6,2){\circle*{0.2}}
  \put(7,3){\circle*{0.2}}
  \put(8,2){\circle*{0.2}}
  \put(9,1){\circle*{0.2}}
  \put(10,0){\circle*{0.2}}
  \put(11,1){\circle*{0.2}}
  \put(12,0){\circle*{0.2}}
  \put(13,1){\circle*{0.2}}
  \put(14,0){\circle*{0.2}}
  \put(15,1){\circle*{0.2}}
  \put(16,0){\circle*{0.2}}
  \put(0,0){\line(1,1){1}}
  \put(1,1){\line(1,-1){1}}
  \put(2,0){\line(1,1){1}}
  \put(3,1){\line(1,-1){1}}
  \put(4,0){\line(1,1){3}}
  \put(7,3){\line(1,-1){3}}
  \put(10,0){\line(1,1){1}}
  \put(11,1){\line(1,-1){1}}
  \put(12,0){\line(1,1){1}}
  \put(13,1){\line(1,-1){1}}
  \put(14,0){\line(1,1){1}}
  \put(15,1){\line(1,-1){1}}
  \put(0,-0.7){\makebox(0,0){$0$}}
  \put(4,-0.7){\makebox(0,0){$2i$}}
  \put(10,-0.7){\makebox(0,0){$2j$}}
  \put(16,-0.7){\makebox(0,0){$2n$}}
 \end{picture}
\end{center}
 In particular, the spectrum of $\; \DD_n \;$ is isomorphic
 to the poset of the intervals of a chain $\; \mathcal{C}_{n-2} \;$ with $\; n-1 \;$ elements,
 i.e. $\; \Spec(\DD_n) \simeq \{ (i,j)\in \mathcal{C}_n^2 \; : \; i \leq j-2 \} \simeq \Int(\mathcal{C}_{n-2})$.
\end{proposition}
\begin{proof}
 A path $\; x \;$ is covered by a path $\; y \;$ if it can be obtained from $\; y \;$
 by changing a peak $\; U D \;$ into a valley $\; D U \,$.
 Dyck paths with a unique pyramid are the only paths for which this operation can be performed just in one way.
 Clearly, every join-irreducible is uniquely determined by the interval corresponding to its pyramids.
\end{proof}

\begin{figure}[h]
\begin{center}
 \setlength{\unitlength}{5mm}
 \begin{picture}(18,4)
 \put(0,0){
 \begin{picture}(16,4)
  \put(0,0){\circle*{0.2}}
  \put(1,1){\circle*{0.2}}
  \put(2,2){\circle*{0.2}}
  \put(3,3){\circle*{0.2}}
  \put(4,2){\circle*{0.2}}
  \put(5,3){\circle*{0.2}}
  \put(6,4){\circle*{0.2}}
  \put(7,3){\circle*{0.2}}
  \put(8,2){\circle*{0.2}}
  \put(10,2){\circle*{0.2}}
  \put(9,1){\circle*{0.2}}
  \put(11,1){\circle*{0.2}}
  \put(12,0){\circle*{0.2}}
  \put(13,1){\circle*{0.2}}
  \put(14,2){\circle*{0.2}}
  \put(15,1){\circle*{0.2}}
  \put(16,2){\circle*{0.2}}
  \put(17,1){\circle*{0.2}}
  \put(18,0){\circle*{0.2}}
  \put(0,0){\line(1,1){3}}
  \put(3,3){\line(1,-1){1}}
  \put(4,2){\line(1,1){2}}
  \put(6,4){\line(1,-1){3}}
  \put(9,1){\line(1,1){1}}
  \put(10,2){\line(1,-1){2}}
  \put(12,0){\line(1,1){2}}
  \put(14,2){\line(1,-1){1}}
  \put(15,1){\line(1,1){1}}
  \put(16,2){\line(1,-1){2}}
  \put(0,0){\dashbox{0.2}(18,0){}}
  \put(1,1){\dashbox{0.2}(10,0){}}
  \put(13,1){\dashbox{0.2}(4,0){}}
  \put(2,2){\dashbox{0.2}(2,0){}}
  \put(4,2){\dashbox{0.2}(4,0){}}
  \put(5,3){\dashbox{0.2}(2,0){}}
 \end{picture}}
 \end{picture}
\end{center}
 \caption{A Dyck path $\; \gamma \;$ of length $\; 18 \;$ with its $\; 9 \;$ tunnels
  (two $0$-tunnels, four $1$-tunnels, two $2$-tunnels and one $3$-tunnel).}
 \label{ExSkewTunnels}
\end{figure}
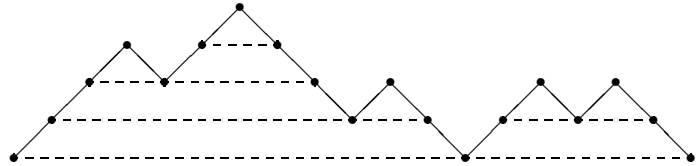

A \emph{$k$-tunnel} of a Dyck path $\; x \;$
is any segment, not reducing to a point, on the horizontal line $\; y = k \;$
having in common with $\; x \;$ only its extreme points \cite{E1,E2}
(see Figure \ref{ExSkewTunnels} for an example).
Clearly every tunnel corresponds to a factor of the form $\; U y D \;$
where the steps $\; U \;$ and $\; D \;$ are at the same level (see again Figure \ref{ExSkewTunnels}).
In particular, the $0$-tunnels correspond to the primitive factors of the paths.

\begin{proposition}\label{meetjoinirredDn}
 The meet of two join-irreducibles of $\; \DD_n \;$ is $\; \widehat{0} \;$ or a join-irreducible,
 i.e. $\; \{ \widehat{0}\} \cup \Spec(\DD_n) \;$ is a ranked sub-meet-semilattice of $\; \DD_n \,$.
 Moreover, the rank of a path in the poset $\; \{ \widehat{0}\} \cup \Spec(\DD_n) \;$
 is equal to the maximum height of its tunnels
 (i.e. $\; r(x) = k \;$ whenever $\; x \;$ has a $k$-tunnel but not a $(k+1)$-tunnel).
\end{proposition}
\begin{proof}
 Two join-irreducible Dyck paths meet in $\; \widehat{0} \;$
 or have pyramids intersecting in a single point with integer coordinates.
\end{proof}

Since Dyck paths are Dyck-like with $\; (a,b) = (1,1) \,$,
from Propositions \ref{young} and \ref{partition} it follows
\begin{proposition}
 The Dyck lattice $\; \DD_n \;$ is isomorphic to the dual of the Young lattice $\; \YY_{\lambda_n} \,$,
 where $\; \lambda_n = (n-1,n-2,\ldots,2,1) \,$.
\end{proposition}

\subsection{Characteristic}\label{chardycklike}

A \emph{Dyck-like lattice} is a distributive lattice
whose spectrum is a ranked poset admitting a labelling of its elements with the following properties:
all labels are positive integers and every antichain $\; S = \{ s_1 ,\ldots, s_n \} \;$ of join-irreducibles
can be linearly ordered so that the labels of the elements of $\; S \;$ are distinct
and, if $\; s_1 \;$ and $\; s_n \;$ are the elements having minimum and maximum labels, respectively,
then $\; s_1 \wedge s_n = s_1 \wedge s_2 \wedge \cdots \wedge s_{n-1} \wedge s_n \,$.
In a Dyck-like lattice, any labelling of the join-irreducibles satisfying the above
properties will be called a \emph{Dyck-like labelling}.
Similarly, if $\; \{ s_1, \ldots, s_n \} \;$ is an antichain of join-irreducibles as above,
the ordered $n$-tuple $\; (s_1, \ldots, s_n) \;$ will be called a \emph{Dyck-like antichain}.

\begin{proposition}
 For any $\; a, b \in \NN \;$ and for every $\; n \in \NN \,$,
 the lattice $\; \DD_n^{(a,b)} \;$ of Dyck-like paths of length $\; n \cdot \ell(a,b) \;$ is a Dyck-like lattice.
\end{proposition}
\begin{proof}
 Consider the lattice $\; \overline{\DD}_n^{(a,b)} \;$ isomorphic to $\; \DD_n^{(a,b)} \,$,
 as defined in Section \ref{sec2}.
 Label each join-irreducible with the abscissa of its unique pyramid.
 Such a labelling is a Dyck-like labelling.
 The fact that the spectrum of $\; \overline{\DD}_n^{(a,b)} \;$ is ranked is a consequence of Theorem \ref{rep}.
\end{proof}

We conjecture that a sort of converse of the previous proposition holds.
More precisely, the following assertion seems plausible:
\begin{conjecture}
 Every finite Dyck-like lattice can be represented
 as a sublattice of a lattice of Dyck-like paths of suitable length.
\end{conjecture}

We will say that an element $\; x \;$ of a distributive lattice $\; \DD \;$ is \emph{quasi-join-irreducible}
when there exists an ordered $k$-tuple $\; (s_1, \ldots, s_k) \;$ forming an antichain of join-irreducibles
such that $\; x = s_1 \vee \cdots \vee s_k \;$ and $\; s_i \wedge s_{i+1} \neq \widehat{0} \,$,
for every $\; i = 1, 2, \ldots, k-1 \,$.

\begin{lemma}\label{Lemma-qji}
 Let $\; \DD \;$ be a Dyck-like lattice and $\; x\in \DD \;$ a quasi-join-irreducible.
 Then $\; x \;$ can be expressed as $\; x = t_1 \vee \cdots \vee t_k \;$
 where $\; (t_1, \ldots, t_k) \;$ is a Dyck-like antichain
 and $\; t_i \wedge t_{i+1} \ne \hat{0} \,$, for every $\; i = 1, 2, \ldots, k -1 \,$.
\end{lemma}
\begin{proof}
 Let $\; x = t_1 \vee \cdots \vee t_k \,$, where $\; (t_1, \ldots, t_k) \;$ is a Dyck-like antichain.
 Suppose there exists an index $\; j < k \;$ such that $\; t_j \wedge t_{j+1} = \hat{0} \,$.
 This would imply $\; t_a\wedge t_b = \hat{0} \,$, for any $\; a \leq j \;$ and $\; b \geq j+1 \,$.
 However, in any rearrangement of the $\; t_i$'s,
 at least a pair of adjacent elements must appear such that one is $\; \leq j \;$ and the other is $\; \geq j+1 \,$,
 and this contradicts the hypothesis that $\; x \;$ is quasi-join-irreducible.
\end{proof}

Lemma \ref{Lemma-qji} asserts that in a Dyck-like lattice
the antichain of join-irreducibles in the definition of a quasi-join-irreducible element
can be taken to be a Dyck-like antichain.

A \emph{special Dyck-like lattice} is a Dyck-like lattice where
the meet of any two join-irreducibles is $\; \widehat{0} \;$ or a join-irreducible.
\begin{proposition}\label{charqji}
 In a special Dyck-like lattice $\; \DD \,$,
 every quasi-join-irreducible element has Euler characteristic equal to $\; 1 \,$.
\end{proposition}
\begin{proof}
 Let $\; x \;$ be a quasi-join-irreducible element of $\; \DD \,$.
 Then $\; x = s_1 \vee \cdots \vee s_k \,$, where $\; s_1 \,$, \ldots, $\; s_k \;$ are incomparable join-irreducibles
 such that $\; s_i \wedge s_{i+1} \neq \widehat{0} \;$ for every $\; i < k \,$.
 If $\; k = 1 \;$ then $\; x \;$ is join-irreducible and $\; \chi(x) = 1 \,$.
 Now we proceed by induction on $\; k \,$. From formula (\ref{dasilva}) we have
 \begin{equation}\label{dasilvaprop}
  \chi(x) = \chi(s_1 \vee \cdots \vee s_k )
          = \sum_{S\subseteq[k]\atop S\neq \emptyset} (-1)^{|S|-1}\chi \left( \bigwedge_{i\in S}s_i \right)\, .
 \end{equation}
 By Lemma \ref{Lemma-qji}, $\; (s_1, \ldots, s_k) \;$ can be taken to be a Dyck-like antichain.
 Hence it follows at once that $\; \bigwedge_{i\in S} s_i = s_{i_1} \wedge s_{i_2} \,$,
 where $\; i_1 = \min S \;$ and $\; i_2 = \max S \,$.
 If $\; S \;$ is a subset of $\; [k] \;$ with $\; \min S = \max S \,$, then clearly $\; |S| = 1 \,$.
 Hence the contribution of these subsets to the sum in (\ref{dasilvaprop}) is $\; \chi(s_1) + \cdots + \chi(s_k) = k \,$.
 If $\; S \;$ is a subset of $\; [k] \;$ with $\; \max S = \min S + 1 \,$, then it follows that $\; S = \{i,i+1\} \,$.
 Since $\; s_i \wedge s_{i+1} \;$ is a join-irreducible, the contribution of these subsets to the sum is
 \begin{displaymath}
  -\sum_{i=1}^{k-1}\chi(s_i \wedge s_{i+1}) = -(k-1) \, .
 \end{displaymath}
 Finally, since the subsets $\; S \;$ of $\; [k] \;$ with $\; \max S - \min S = j \geq 2 \,$,
 having minimum $\; i \,$, maximum $\; i + j \;$ and cardinality $\; h + 2 \;$ ($\, h \geq 0 \,$)
 are exactly $\; { j - 1 \choose h } \,$, the contribution of all these subsets is
 \begin{displaymath}
  \sum_{h=0}^{j-1} { j - 1 \choose h } (-1)^{h+1} \chi(s_i \wedge s_{i+j}) = -(1-1)^{j-1} \chi(s_i \wedge s_{i+j}) = 0 \, .
 \end{displaymath}
 In conclusion, we have $\; \chi(x) = k - ( k - 1 ) = 1 \,$.
\end{proof}

An element $\; x \;$ of a finite distributive lattice is said to have a \emph{quasi-join-irreducible decomposition}
when it can be expressed as a join of quasi-join-irreducible elements
$\; x_1 \,$, \ldots, $\; x_k \,$ such that $\; x_i \wedge x_j = \widehat{0} \,$, for every $\; i \neq j \,$.
\begin{proposition}\label{qjid}
 Every element of a Dyck-like lattice has a quasi-join-irreducible decomposition.
\end{proposition}
\begin{proof}
 Let $\; x = s_1 \vee \cdots \vee s_k \;$ be a join-irreducible decomposition of $\; x \,$,
 where $\; (s_1, \ldots, s_k) \;$ is a Dyck-like antichain.
 If $\; i \;$ is the first index such that $\; s_i \wedge s_{i+1} = \widehat{0} \;$
 then $\; s_1 \wedge \cdots \wedge s_i \;$ is quasi-join-irreducible.
 If $\; j \;$ is the second index satisfying the above condition (and so $\; s_j \wedge s_{j+1} = \widehat{0} \,$),
 then obviously $\; s_{i+1} \wedge \cdots \wedge s_j \;$ is quasi-join-irreducible
 and $\; ( s_1 \wedge \cdots \wedge s_i ) \wedge ( s_{i+1} \wedge \cdots \wedge s_j )= \widehat{0} \,$.
 Repeating this argument one obtains the desired decomposition.
\end{proof}

\begin{theorem}\label{char}
 Let $\; \DD \;$ be a special Dyck-like lattice.
 Then, for every $\; x \in \DD \,$, $\; \chi(x) \;$ is the number of quasi-join-irreducibles in a decomposition of $\; x \,$.
\end{theorem}
\begin{proof}
 If $\; x = x_1 \vee \cdots \vee x_k \;$ is a quasi-join-irreducible decomposition of $\; x \,$,
 then it follows at once that $\; \chi(x) = \chi(x_1) + \cdots + \chi(x_k) = k \,$.
\end{proof}

\begin{corollary}
 Let $\; \DD \;$ be a special Dyck-like lattice.
 Then two quasi-join-irreducible decompositions of the same element $\; x \;$ have the same number of elements.
\end{corollary}

In the case of Dyck paths, the characteristic can be interpreted combinatorially as follows.
\begin{proposition}\label{dyckqji}
 A Dyck path $\; x \in \DD_n \;$ is quasi-join-irreducible
 if and only if it has precisely one $1$-tunnel.
\end{proposition}
\begin{proof}
 If $\; x \;$ has precisely one $1$-tunnel,
 then it is of the form $\; (UD)^s y(UD)^t \,$, for suitable $\; s, t \in \NN \,$,
 where $\; y \;$ is an elevated Dyck path of length $\; > 2 \,$.
 Then $\; x \;$ can be expressed as $\; x = x_1 \vee \cdots \vee x_k \,$,
 where $\; x_1 \,$, \ldots, $\; x_k \;$ are the join-irreducibles uniquely determined by the peaks of $\; y \,$.
 The fact that $\; y \;$ is elevated implies that $\; x_i \wedge x_{i+1} \neq \widehat{0} \,$, for every $\; i < k \,$,
 and consequently that $\; x \;$ is quasi-join-irreducible.

 On the other hand, suppose that $\; x \in \DD_n \;$ is quasi-join-irreducible.
 If $\; x \;$ had no $1$-tunnels, then $\; x = \widehat{0} \,$,
 which is impossible (since $\; \widehat{0} \;$ is not quasi-join-irreducible).
 If $\; x \;$ had more than one $1$-tunnel, then $\; x \;$ would have at least two elevated factors,
 that is $\; x = \alpha v \beta w \gamma \,$, with $\; v \;$ and $\; w \;$ elevated Dyck paths.
 In this situation, any expression of $\; x \;$ as a join of join-irreducibles
 would contain join-irreducible elements determined by the peaks of all the elevated factors of $\; x \,$.
 Thus, if $\; x = x_1 \vee \cdots \vee x_k \;$ is any join-irreducible decomposition of $\; x \,$,
 then there exists at least one $\; i < k \;$ such that $\; x_i \;$ and $\; x_{i+1} \;$ are join-irreducibles
 determined by the peaks of two different elevated factors,
 and so $\; x_i \wedge x_{i+1} = \widehat{0} \,$.
 In conclusion, $\; x \;$ has exactly one $1$-tunnel.
\end{proof}

Finally, as an immediate consequence of Theorem \ref{char} and Proposition \ref{dyckqji}, we have
\begin{theorem}\label{chardyck}
 The characteristic of a Dyck path is the number of its $1$-tunnels.
\end{theorem}

\subsection{Generalized characteristics}

Suppose that $\; \DD \;$ is a finite special Dyck-like lattice with spectrum $\; P \;$
such that $\; \widehat{P} = \{ \widehat{0} \} \cup P \;$ is ranked with rank function $\; r_{\widehat{P}} \,$.
For every $\; k \in \NN \,$, we define the \emph{generalized characteristic} $\; \chi_k \;$
as the valuation on $\; \DD \;$ such that
\begin{displaymath}
 \chi_k (x)=
  \begin{cases}
   1 & \textrm{if}\;\; r_{\widehat{P}}(x)\geq k\\
   0 & \textrm{if}\;\; r_{\widehat{P}}(x)<k
  \end{cases}
\end{displaymath}
for every join-irreducible $\; x \,$.
Clearly $\; \chi_1 \;$ is equal to the ordinary characteristic $\; \chi \,$.
Our aim is to evaluate $\; \chi_k (x) \;$ for every $\; x \in \DD \,$.
\begin{proposition}\label{propsemijoinirredk}
 Let $\; x = x_1 \vee \cdots \vee x_m \;$ where $\; m \geq 1 \,$, each $\; x_i \;$ is a join-irreducible
 and $\; r_{\widehat{P}}(x_j \wedge x_{j+1}) \geq k \;$ for every $\; j < m \,$.
 Then $\; \chi_k(x) = 1 \,$.
\end{proposition}
\begin{proof}
 The proof follows the same lines of that of Proposition \ref{charqji}.
\end{proof}

\begin{proposition}\label{chikjoinsemiirreducibles}
 Let $\; x = x_1 \vee \cdots \vee x_m \,$, where $\; m \geq 1 \,$,
 and each $\; x_i = x_{i,1} \vee \cdots \vee x_{i,\ell_i} \;$ is a join of join-irreducibles
 such that $\; r_{\widehat{P}}(x_{i,j} \wedge x_{i,j+1}) \geq k \,$, for every $\; j = 1, \ldots, \ell_i-1 \,$,
 and $\; r_{\widehat{P}}(x_i\wedge x_j ) < k \;$ for every $\; i \neq j \,$.
 Then $\; \chi_k(x) = m \,$.
\end{proposition}
\begin{proof}
 Using formula (\ref{dasilva}), we have
 \begin{displaymath}
  \chi_k(x) = \chi_k(x_1 \vee \cdots \vee x_m )
  = \sum_{S\subseteq [m]\atop S\neq \emptyset}(-1)^{|S|-1}\chi \left( \bigwedge_{i\in S}x_i \right)\, .
 \end{displaymath}
 By hypothesis $\; \bigwedge_{i\in S}x_i \;$ is a join-irreducible
 with $\; r_{\widehat{P}}\left( \bigwedge_{i\in S}x_i \right) < k \;$
 for every $\; S \subseteq [m] \,$, $\; |S| \geq 2 \,$.
 Hence in these cases $\; \chi(\bigwedge_{i\in S}x_i ) = 0 \;$
 and then $\; \chi_k(x) = \chi_k(x_1) + \cdots + \chi_k(x_m) \,$.
 Finally, the claim follows applying Proposition \ref{propsemijoinirredk}.
\end{proof}

\begin{proposition}\label{kjoinsemiirreducibledecomposition}
 Every element $\; x \;$ of $\; \DD \;$ can be written as
 $\; x = ( x_1 \vee \cdots \vee x_h ) \vee ( x_{h+1} \vee \cdots \vee x_m ) \;$ where
 \begin{enumerate}
  \item
   $ x_i = x_{i,1} \vee \cdots \vee x_{i,\ell_i} \;$
   is a join of join-irreducibles with $\; r_{\widehat{P}}(x_{i,j} \wedge x_{i,j+1}) \geq k \;$
   for every $\; j = 1, \ldots, \ell_i-1 \,$,
   and $\; r_{\widehat{P}}( x_i \wedge x_j ) < k \;$ whenever $\; i \neq j \,$,
   for every $\; i = 1,\ldots, h \,$;
  \item
   $\; x_i \;$ is a join-irreducible with $\; r_{\widehat{P}}(x_i) < k \,$, for every $\; i = h + 1, \ldots, m \,$.
 \end{enumerate}
\end{proposition}
\begin{proof}
 Let $\; x = p_1 \vee \cdots \vee p_s \,$, where $\; (p_1, \ldots, p_s) \;$ is a Dyck-like antichain.
 For the first element $\; p_1 \;$ there are two possible cases.
 If $\; r_{\widehat{P}}(p_1) < k \,$, then $\; p_1 \;$ is one of the $\; x_i$'s.
 Otherwise, if $\; r_{\widehat{P}}(p_1) \geq k \,$,
 then consider the first index $\; i \;$ such that $\; r_{\widehat{P}}(p_i \wedge p_{i+1}) < k \,$:
 then $\; p_1 \vee \cdots \vee p_i \;$ is a join of join-irreducibles such that
 $\; r_{\widehat{P}}(p_j \wedge p_{j+1}) \geq k \,$, for $\; j < i \,$,
 and so it is one of the $\; x_i$'s.
 Repeating this argument and rearranging the $\; x_i \;$ in the correct order,
 we obtain the desired decomposition.
\end{proof}

Any decomposition of the kind described in Proposition \ref{kjoinsemiirreducibledecomposition}
will be said a \emph{$k$-quasi-join-irreducible decomposition} of the element $\; x \in \DD \,$,
and the elements $\; x_1 \,$, \ldots, $\; x_h \;$ appearing in such a decomposition
will be called \emph{$k$-quasi-join-irreducibles}.

\begin{theorem}\label{chikthm}
 The generalized characteristic $\; \chi_k(x) \;$ of an element $\; x \in \DD \;$
 is equal to the number of $k$-quasi-join-irreducibles in any $k$-quasi-join-irreducible decomposition of $\; x \,$.
\end{theorem}
\begin{proof}
 By Proposition \ref{kjoinsemiirreducibledecomposition},
 every element $\; x \;$ of $\; \DD \;$ admits a $k$-quasi-join-irreducible decomposition
 $\; x = x_1 \vee \cdots \vee x_h \vee x_{h+1}\vee \cdots \vee x_m \,$.
 Applying formula (\ref{dasilva}) we have
 \begin{displaymath}
  \chi_k(x) = \sum_{S\subseteq [m]\atop S\neq \emptyset } (-1)^{|S|-1}\chi_k \left( \bigwedge_{i\in S}x_i \right)\, .
 \end{displaymath}
 If $\; S \;$ contains an $\; i \;$ such that $\; r_{\widehat{P}}(x_i) < k \;$
 then clearly $\; r_{\widehat{P}}\left( \bigwedge_{i\in S} x_i \right) < k \;$
 and $\; \chi_k\left( \bigwedge_{i\in S} x_i \right) = 0 \,$.
 Therefore in the computation of $\; \chi_k(x) \;$ all the join-irreducibles
 with rank strictly less than $\; k \;$ in $\; \widehat{P} \;$ can be discarded, i.e.
 \begin{displaymath}
  \chi_k(x) = \sum_{S\subseteq [h]\atop S\neq \emptyset} (-1)^{|S|-1}\chi_k \left( \bigwedge_{i\in S}x_i \right)
            = \chi_k(x_1 \vee \cdots \vee x_h ).
 \end{displaymath}
 Finally, the claim follows from Proposition \ref{chikjoinsemiirreducibles}.
\end{proof}

As a consequence of Theorem \ref{chikthm} it follows that the number of $k$-quasi-join-irreducibles
in any $k$-quasi-join-irreducible decomposition of an element $\; x \;$ is constant.
Moreover, from Theorem \ref{chikthm} and Proposition \ref{meetjoinirredDn},
we have the following interpretation of the generalized characteristics of Dyck lattices.
\begin{theorem}
 The generalized characteristic $\; \chi_k(x) \;$ of an element $\; x \;$ of the Dyck lattice $\; \DD_n \;$
 is equal to the number of $k$-tunnels of $\; x \,$.
\end{theorem}

\section{Motzkin lattices}

\subsection{Representation}\label{motzkin}

\begin{figure}[h]
\begin{center}
 \setlength{\unitlength}{5mm}
 \begin{picture}(28,10)
  \put(0,6){
   \begin{picture}(3,4)
    \put(1,0){\circle*{0.2}}
    \put(0,1){\circle*{0.2}} \put(0,1){\circle{0.4}}
    \put(1,1){\circle*{0.2}} \put(1,1){\circle{0.4}}
    \put(2,1){\circle*{0.2}} \put(2,1){\circle{0.4}}
    \put(0,2){\circle*{0.2}}
    \put(1,2){\circle*{0.2}}
    \put(2,2){\circle*{0.2}}
    \put(1,3){\circle*{0.2}}
    \put(1,4){\circle*{0.2}} \put(1,4){\circle{0.4}}
    \put(1,0){\line(-1,1){1}}
    \put(1,0){\line(0,1){1}}
    \put(1,0){\line(1,1){1}}
    \put(0,1){\line(0,1){1}}
    \put(0,1){\line(1,1){1}}
    \put(1,1){\line(-1,1){1}}
    \put(1,1){\line(1,1){1}}
    \put(2,1){\line(-1,1){1}}
    \put(2,1){\line(0,1){1}}
    \put(1,3){\line(-1,-1){1}}
    \put(1,3){\line(0,-1){1}}
    \put(1,3){\line(1,-1){1}}
    \put(1,3){\line(0,1){1}}
    \put(2.2,3.1){\makebox(0,0){\small$\MM_4$}}
   \end{picture}}
  \put(2,5.7){
   \begin{picture}(3,2)
    \put(0,0){\circle*{0.2}}
    \put(1,0){\circle*{0.2}}
    \put(2,0){\circle*{0.2}}
    \put(1,1){\circle*{0.2}}
    \put(1,1){\line(-1,-1){1}}
    \put(1,1){\line(0,-1){1}}
    \put(1,1){\line(1,-1){1}}
    \put(1,-0.8){\makebox(0,0){\small$\Spec(\MM_4)$}}
   \end{picture}}
  \put(0,0){
   \begin{picture}(3,4)
    \put(1,0){\circle*{0.2}} \put(1,-0.5){\makebox(0,0){\footnotesize$0$}}
    \put(0,1){\circle*{0.2}} \put(-0.4,1){\makebox(0,0){\footnotesize$1$}}
    \put(1,1){\circle*{0.2}} \put(1.4,1){\makebox(0,0){\footnotesize$1$}}
    \put(2,1){\circle*{0.2}} \put(2.4,1){\makebox(0,0){\footnotesize$1$}}
    \put(0,2){\circle*{0.2}} \put(-0.4,2){\makebox(0,0){\footnotesize$2$}}
    \put(1,2){\circle*{0.2}} \put(1.4,2){\makebox(0,0){\footnotesize$2$}}
    \put(2,2){\circle*{0.2}} \put(2.4,2){\makebox(0,0){\footnotesize$2$}}
    \put(1,3){\circle*{0.2}} \put(1.4,3){\makebox(0,0){\footnotesize$3$}}
    \put(1,4){\circle*{0.2}} \put(1.4,4){\makebox(0,0){\footnotesize$1$}}
    \put(1,0){\line(-1,1){1}}
    \put(1,0){\line(0,1){1}}
    \put(1,0){\line(1,1){1}}
    \put(0,1){\line(0,1){1}}
    \put(0,1){\line(1,1){1}}
    \put(1,1){\line(-1,1){1}}
    \put(1,1){\line(1,1){1}}
    \put(2,1){\line(-1,1){1}}
    \put(2,1){\line(0,1){1}}
    \put(1,3){\line(-1,-1){1}}
    \put(1,3){\line(0,-1){1}}
    \put(1,3){\line(1,-1){1}}
    \put(1,3){\line(0,1){1}}
   \end{picture}}
  \put(5,0){
   \begin{picture}(11,10)
    \put( 5,0){\circle*{0.2}}
    \put( 1,2){\circle*{0.2}} \put(1,2){\circle{0.4}}
    \put( 3,2){\circle*{0.2}} \put(3,2){\circle{0.4}}
    \put( 7,2){\circle*{0.2}} \put(7,2){\circle{0.4}}
    \put( 9,2){\circle*{0.2}} \put(9,2){\circle{0.4}}
    \put( 0,4){\circle*{0.2}}
    \put( 2,4){\circle*{0.2}}
    \put( 4,4){\circle*{0.2}}
    \put( 6,4){\circle*{0.2}}
    \put( 8,4){\circle*{0.2}}
    \put(10,4){\circle*{0.2}}
    \put( 1,6){\circle*{0.2}}
    \put( 3,6){\circle*{0.2}}
    \put( 7,6){\circle*{0.2}}
    \put( 9,6){\circle*{0.2}}
    \put( 0,8){\circle*{0.2}} \put(0,8){\circle{0.4}}
    \put( 5,8){\circle*{0.2}}
    \put(10,8){\circle*{0.2}} \put(10,8){\circle{0.4}}
    \put( 3,9){\circle*{0.2}}
    \put( 7,9){\circle*{0.2}}
    \put(5,10){\circle*{0.2}}
    \put(5,0){\line(-2,1){4}}
    \put(5,0){\line(-1,1){2}}
    \put(5,0){\line(1,1){2}}
    \put(5,0){\line(2,1){4}}
    \put(1,2){\line(-1,2){1}}
    \put(1,2){\line(1,2){1}}
    \put(1,2){\line(3,2){3}}
    \put(3,2){\line(-3,2){3}}
    \put(3,2){\line(3,2){3}}
    \put(3,2){\line(5,2){5}}
    \put(7,2){\line(-5,2){5}}
    \put(7,2){\line(-1,2){1}}
    \put(7,2){\line(3,2){3}}
    \put(9,2){\line(-5,2){5}}
    \put(9,2){\line(-1,2){1}}
    \put(9,2){\line(1,2){1}}
    \put(0,4){\line(1,2){1}}
    \put(0,4){\line(3,2){3}}
    \put(2,4){\line(-1,2){1}}
    \put(2,4){\line(5,2){5}}
    \put(4,4){\line(-1,2){1}}
    \put(4,4){\line(3,2){3}}
    \put(6,4){\line(-5,2){5}}
    \put(6,4){\line(3,2){3}}
    \put(8,4){\line(-5,2){5}}
    \put(8,4){\line(1,2){1}}
    \put(10,4){\line(-3,2){3}}
    \put(10,4){\line(-1,2){1}}
    \put(1,6){\line(-1,2){1}}
    \put(1,6){\line(2,1){4}}
    \put(3,6){\line(1,1){2}}
    \put(7,6){\line(-1,1){2}}
    \put(9,6){\line(-2,1){4}}
    \put(9,6){\line(1,2){1}}
    \put(0,8){\line(3,1){3}}
    \put(5,8){\line(-2,1){2}}
    \put(5,8){\line(2,1){2}}
    \put(10,8){\line(-3,1){3}}
    \put(3,9){\line(2,1){2}}
    \put(7,9){\line(-2,1){2}}
    \put(2,0.5){\makebox(0,0){\small$\MM_5$}}
   \end{picture}}
  \put(14.5,0.5){
   \begin{picture}(3,2)
    \put(0,0){\circle*{0.2}}
    \put(1,0){\circle*{0.2}}
    \put(2,0){\circle*{0.2}}
    \put(3,0){\circle*{0.2}}
    \put(1,1){\circle*{0.2}}
    \put(2,1){\circle*{0.2}}
    \put(1,1){\line(-1,-1){1}}
    \put(1,1){\line(0,-1){1}}
    \put(1,1){\line(1,-1){1}}
    \put(2,1){\line(-1,-1){1}}
    \put(2,1){\line(0,-1){1}}
    \put(2,1){\line(1,-1){1}}
    \put(1.5,-0.8){\makebox(0,0){\small$\Spec(\MM_5)$}}
   \end{picture}}
  \put(17,0){
   \begin{picture}(11,10)
    \put( 5,0){\circle*{0.2}} \put(5,-0.5){\makebox(0,0){\footnotesize$0$}}
    \put( 1,2){\circle*{0.2}} \put(0.5,2){\makebox(0,0){\footnotesize$1$}}
    \put( 3,2){\circle*{0.2}} \put(2.5,2){\makebox(0,0){\footnotesize$1$}}
    \put( 7,2){\circle*{0.2}} \put(7.5,2){\makebox(0,0){\footnotesize$1$}}
    \put( 9,2){\circle*{0.2}} \put(9.5,2){\makebox(0,0){\footnotesize$1$}}
    \put( 0,4){\circle*{0.2}} \put(-0.5,4){\makebox(0,0){\footnotesize$2$}}
    \put( 2,4){\circle*{0.2}} \put(1.5,4){\makebox(0,0){\footnotesize$2$}}
    \put( 4,4){\circle*{0.2}} \put(3.5,4){\makebox(0,0){\footnotesize$2$}}
    \put( 6,4){\circle*{0.2}} \put(6.5,4){\makebox(0,0){\footnotesize$2$}}
    \put( 8,4){\circle*{0.2}} \put(8.5,4){\makebox(0,0){\footnotesize$2$}}
    \put(10,4){\circle*{0.2}} \put(10.5,4){\makebox(0,0){\footnotesize$2$}}
    \put( 1,6){\circle*{0.2}} \put(0.5,6){\makebox(0,0){\footnotesize$3$}}
    \put( 3,6){\circle*{0.2}} \put(2.5,6){\makebox(0,0){\footnotesize$3$}}
    \put( 7,6){\circle*{0.2}} \put(7.5,6){\makebox(0,0){\footnotesize$3$}}
    \put( 9,6){\circle*{0.2}} \put(9.5,6){\makebox(0,0){\footnotesize$3$}}
    \put( 0,8){\circle*{0.2}} \put(-0.5,8){\makebox(0,0){\footnotesize$1$}}
    \put( 5,8){\circle*{0.2}} \put(5,8.5){\makebox(0,0){\footnotesize$4$}}
    \put(10,8){\circle*{0.2}} \put(10.5,8){\makebox(0,0){\footnotesize$1$}}
    \put( 3,9){\circle*{0.2}} \put(3,9.5){\makebox(0,0){\footnotesize$2$}}
    \put( 7,9){\circle*{0.2}} \put(7,9.5){\makebox(0,0){\footnotesize$2$}}
    \put(5,10){\circle*{0.2}} \put(5,10.5){\makebox(0,0){\footnotesize$0$}}
    \put(5,0){\line(-2,1){4}}
    \put(5,0){\line(-1,1){2}}
    \put(5,0){\line(1,1){2}}
    \put(5,0){\line(2,1){4}}
    \put(1,2){\line(-1,2){1}}
    \put(1,2){\line(1,2){1}}
    \put(1,2){\line(3,2){3}}
    \put(3,2){\line(-3,2){3}}
    \put(3,2){\line(3,2){3}}
    \put(3,2){\line(5,2){5}}
    \put(7,2){\line(-5,2){5}}
    \put(7,2){\line(-1,2){1}}
    \put(7,2){\line(3,2){3}}
    \put(9,2){\line(-5,2){5}}
    \put(9,2){\line(-1,2){1}}
    \put(9,2){\line(1,2){1}}
    \put(0,4){\line(1,2){1}}
    \put(0,4){\line(3,2){3}}
    \put(2,4){\line(-1,2){1}}
    \put(2,4){\line(5,2){5}}
    \put(4,4){\line(-1,2){1}}
    \put(4,4){\line(3,2){3}}
    \put(6,4){\line(-5,2){5}}
    \put(6,4){\line(3,2){3}}
    \put(8,4){\line(-5,2){5}}
    \put(8,4){\line(1,2){1}}
    \put(10,4){\line(-3,2){3}}
    \put(10,4){\line(-1,2){1}}
    \put(1,6){\line(-1,2){1}}
    \put(1,6){\line(2,1){4}}
    \put(3,6){\line(1,1){2}}
    \put(7,6){\line(-1,1){2}}
    \put(9,6){\line(-2,1){4}}
    \put(9,6){\line(1,2){1}}
    \put(0,8){\line(3,1){3}}
    \put(5,8){\line(-2,1){2}}
    \put(5,8){\line(2,1){2}}
    \put(10,8){\line(-3,1){3}}
    \put(3,9){\line(2,1){2}}
    \put(7,9){\line(-2,1){2}}
   \end{picture}}
 \end{picture}
\end{center}
\caption{The Motzkin lattices $\; \MM_4 \;$ and $\; \MM_5 \,$, their spectra and the distribution of the characteristic.}
\label{MotzkinLatticeM5}
\end{figure}
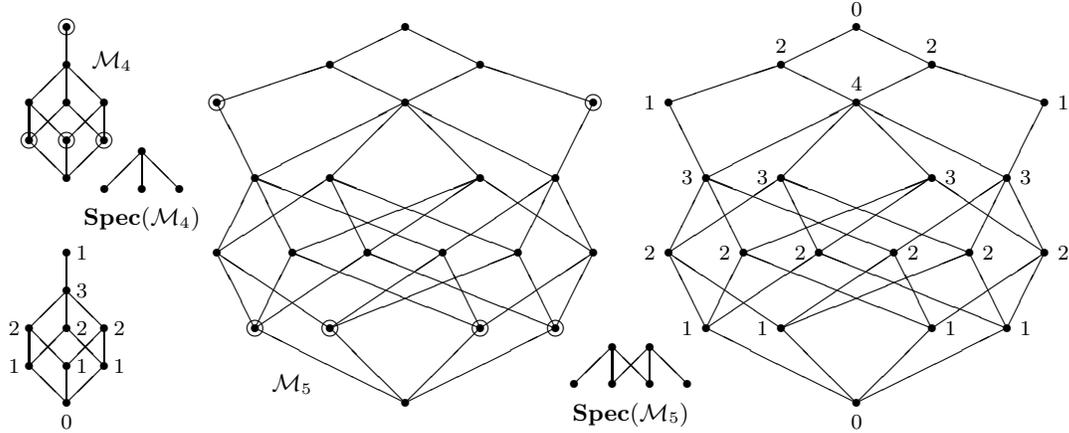

The \emph{Motzkin lattice} $\; \MM_n \;$ is a distributive lattice
of height $\; \left\lfloor\frac{n}{2}\right\rfloor \left\lceil\frac{n}{2}\right\rceil \,$, with minimum $\; H^n \,$,
and maximum $\; U^k D^k \;$ when $\; n = 2k \;$ or $\; U^k H D^k \;$ when $\; n = 2 k + 1 \,$,
with $\; n - 1 \;$ atoms of the form $\; H^k U D H^{n-k-2} \;$
and just one coatom $\; U^{k-1} H D^{k-1} \;$ when $\; n = 2k $
and two coatoms $\; U^{k-1} H U D^k \;$ and $\; U^k D H D^{k-1} \;$ when $\; n = 2 k + 1 \,$,
It is easy to see that the socle is the path $\; s = UH^{n-2}D \;$
and that the principal ideal $\; \downarrow s \;$ is isomorphic to a Boolean algebra $\; B_{n-1} \;$
while the principal filter $\; \uparrow s \;$ is isomorphic to a Motzkin lattice $\; \MM_{n-2} \,$,
whenever $\; n \geq 1 \;$ (see Figure \ref{MotzkinLatticeM5}).
The rank function is given by the area determined by the path, i.e. $\; r(x) = \AA(x) \,$.

\begin{proposition}
 The join-irreducibles of the Motzkin lattice $\; \MM_n \;$ are the paths with exactly one peak, i.e. having the form
 \begin{center}
 \setlength{\unitlength}{3mm}
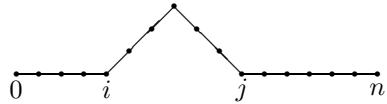
 \begin{picture}(16,4)(0,-1)
  \put(0,0){\circle*{0.2}}
  \put(1,0){\circle*{0.2}}
  \put(2,0){\circle*{0.2}}
  \put(3,0){\circle*{0.2}}
  \put(4,0){\circle*{0.2}}
  \put(5,1){\circle*{0.2}}
  \put(6,2){\circle*{0.2}}
  \put(7,3){\circle*{0.2}}
  \put(8,2){\circle*{0.2}}
  \put(9,1){\circle*{0.2}}
  \put(10,0){\circle*{0.2}}
  \put(11,0){\circle*{0.2}}
  \put(12,0){\circle*{0.2}}
  \put(13,0){\circle*{0.2}}
  \put(14,0){\circle*{0.2}}
  \put(15,0){\circle*{0.2}}
  \put(16,0){\circle*{0.2}}
  \put(0,0){\line(1,0){4}}
  \put(10,0){\line(1,0){6}}
  \put(4,0){\line(1,1){3}}
  \put(7,3){\line(1,-1){3}}
  \put(0,-0.7){\makebox(0,0){$0$}}
  \put(4,-0.7){\makebox(0,0){$i$}}
  \put(10,-0.7){\makebox(0,0){$j$}}
  \put(16,-0.7){\makebox(0,0){$n$}}
 \end{picture}
\end{center}
 The spectrum of $\; \MM_n \;$ is isomorphic to the poset of the intervals of even length
 of a chain having $\; n + 1 \;$ elements, i.e.
 $\; \Spec(\MM_n) \simeq \{ (i,j) \in \CC_n^2 \; : \; \exists k \in \NN \; ( j - i = 2 k + 2 ) \}$.
\end{proposition}
\begin{proof}
 If $\; x \in \MM_n \;$ has an horizontal step at height $\; > 0 \,$,
 then it is easy to see that $\; x \;$ can be obtained as the join of two smaller paths.
 Hence a join-irreducible cannot have horizontal steps at height $\; > 0 \,$.
 If $\; x \;$ had more than one peak, then it could be expressed as the join of all paths obtained from $\; x \;$
 by replacing each peak with a couple of horizontal steps, one peak at a time.
 So $\; x \;$ can have only one peak.
\end{proof}

\begin{remark}
 Motzkin paths are not Dyck-like paths,
 nevertheless the lattice $\; \MM_n \;$ of Motzkin paths of length $\; n \;$ is a Dyck-like lattice.
 Notice that $\; \MM_n \;$ is isomorphic to the lattice of Dyck paths of length $\; 2 n \;$
 having at most two consecutive down steps \cite{BF}.
 This agrees with our previous conjecture on the representation of Dyck-like lattices.

 Motzkin lattices are not special Dyck-like lattices.
 For instance, in $\; \MM_5 \;$ the meet of the paths $\; U^2 D^2 H \;$ and $\; H U^2 D^2 \;$ is $\; HUHDH \,$,
 which is different from the minimum $\; \widehat{0} \;$ and non join-irreducible.
 However, we can prove a result similar to Proposition \ref{meetjoinirredDn}
 for the meet of two join-irreducibles in $\; \MM_n \,$,
 which will allow to compute the Euler characteristic also in this case.
\end{remark}

\subsection{Characteristic}

To give a combinatorial description of the characteristic for Motzkin lattices
we cannot use the theory developed in the previous sections,
since in a Motzkin lattice it could happen that the meet of two join-irreducibles
is neither $\; \widehat{0} \;$ nor a join-irreducible, as we have seen in section \ref{motzkin}.
However, the arguments developed for Dyck paths can be adapted to the Motzkin case and lead to analogous results.

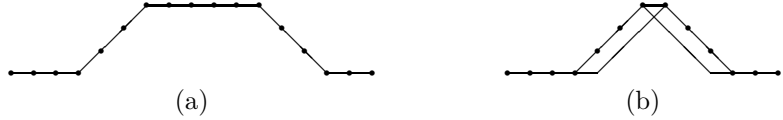
\begin{figure}[h]
 \begin{center}
 \setlength{\unitlength}{3mm}
 \begin{picture}(34,3)
 \put(0,0){
 \begin{picture}(16,3)
  \put(0,0){\circle*{0.2}}
  \put(1,0){\circle*{0.2}}
  \put(2,0){\circle*{0.2}}
  \put(3,0){\circle*{0.2}}
  \put(4,1){\circle*{0.2}}
  \put(5,2){\circle*{0.2}}
  \put(6,3){\circle*{0.2}}
  \put(7,3){\circle*{0.2}}
  \put(8,3){\circle*{0.2}}
  \put(9,3){\circle*{0.2}}
  \put(10,3){\circle*{0.2}}
  \put(11,3){\circle*{0.2}}
  \put(12,2){\circle*{0.2}}
  \put(13,1){\circle*{0.2}}
  \put(14,0){\circle*{0.2}}
  \put(15,0){\circle*{0.2}}
  \put(16,0){\circle*{0.2}}
  \put(0,0){\line(1,0){3}}
  \put(3,0){\line(1,1){3}}
  \put(6,3){\line(1,0){5}}
  \put(11,3){\line(1,-1){3}}
  \put(14,0){\line(1,0){2}}
  \put(8,-1.3){\makebox(0,0){(a)}}
 \end{picture}}
 \put(22,0){
 \begin{picture}(12,3)
  \put(0,0){\circle*{0.2}}
  \put(1,0){\circle*{0.2}}
  \put(2,0){\circle*{0.2}}
  \put(3,0){\circle*{0.2}}
  \put(4,1){\circle*{0.2}}
  \put(5,2){\circle*{0.2}}
  \put(6,3){\circle*{0.2}}
  \put(7,3){\circle*{0.2}}
  \put(8,2){\circle*{0.2}}
  \put(9,1){\circle*{0.2}}
  \put(10,0){\circle*{0.2}}
  \put(11,0){\circle*{0.2}}
  \put(12,0){\circle*{0.2}}
  \put(0,0){\line(1,0){3}}
  \put(3,0){\line(1,1){3}}
  \put(6,3){\line(1,0){1}}
  \put(7,3){\line(1,-1){3}}
  \put(10,0){\line(1,0){2}}
  \put(6,3){\line(1,-1){3}}
  \put(7,3){\line(-1,-1){3}}
  \put(3,0){\line(1,0){1}}
  \put(9,0){\line(1,0){1}}
  \put(6,-1.3){\makebox(0,0){(b)}}
 \end{picture}}
 \end{picture}
\end{center}
\caption{(a) A Motzkin path in $\; T_{16,5,3} \,$.
(b) A Motzkin path in $\; T_{12,1,3} \;$ and the two join-irreducibles it covers.}
\label{TruncatedPyramidEx}
\end{figure}

A \emph{truncated pyramid} of a Motzkin path $\; x \;$
is a sequence of $\; k \geq 1 \;$ up steps followed by a sequence of $\; m \geq 1 \;$ horizontal steps
followed by a sequence of $\; k \;$ down steps, i.e. $\; U^k H^m D^k \,$.
The positive integer $\; k \;$ is called the \emph{dimension} of the truncated pyramid,
whereas $\; m \;$ is its \emph{length}.
Moreover, we say that a truncated pyramid has \emph{height} $\; h \;$
when the sequence of horizontal steps lies on the line $\; y = h \,$.
In the sequel, we will denote by $\; T_{n,m,k} \;$ the set of Motzkin paths of length $\; n \;$
having only horizontal steps at height $\; 0 \,$,
except for a unique truncated pyramid of dimension $\; k \;$ and length $\; m \,$.
An element of $\; T_{n,m,k} \;$ will be called
a \emph{Motzkin path with a unique truncated pyramid} of length $\; m \;$ and dimension $\; k \;$
(see Figure \ref{TruncatedPyramidEx}).

Now we are ready to state our result on the meet of two join-irreducible Motzkin paths.
\begin{proposition}
 In $\; \MM_n \;$ the meet of two join-irreducibles is either a join-irreducible or an element of $\; T_{n,1,k} \;$
 (i.e., a Motzkin path with a unique truncated pyramid of length $\; 1 \,$).
\end{proposition}
\begin{proof}
 Let $\; E_n \;$ be the set obtained by taking all join-irreducibles of $\; \MM_n \;$
 and all elements of $\; T_{n,1,h} \,$, for $\; h = 1, 2, \ldots, n - 2 \,$.
 $\; E_n \;$ with the induced order is isomorphic to $\; \Spec(\DD_n) \,$.
 Indeed, in $\; E_n \;$ each element of $\; T_{n,1,h} \;$ covers precisely two join-irreducibles
 (see Figure \ref{TruncatedPyramidEx}(b))
 and each join-irreducible covers precisely two elements of $\; T_{n,1,h} \,$, for a suitable $\; h \,$.
 Alternatively, map each element of $\; \Spec(\DD_n) \;$ to the element of $\; E_n \;$
 obtained by performing the following substitutions, when reading the path from left to right:
 $\; UU \rightarrow U \,$, $\; UD \rightarrow H \,$, $\; DD \rightarrow D \,$.
 The resulting map is an isomorphism. Hence the claim follows from Proposition \ref{meetjoinirredDn}.
\end{proof}

First of all, we compute the characteristic of some particular Motzkin paths.
\begin{lemma}\label{truncpyr}
 If $\; x \in T_{n,m,h} \,$, then $\; \chi(x) = (-1)^{h+1} m + 1 \,$.
\end{lemma}
\begin{proof}
 Since a Motzkin path $\; x \in T_{n,1,1} \;$ is the join of two join-irreducibles whose meet is $\; \widehat{0} \,$,
 we have $\; \chi(x) = 2 \,$.
 Similarly, since a Motzkin path $\; x \in T_{n,1,2} \;$
 is the join of two join-irreducibles whose meet is a Motzkin path belonging to $\; T_{n,1,1} \,$,
 we have $\; \chi(x) = 1 + 1 - 2 = 0 \,$.
 Iterating this argument it follows that for every $\; x \in T_{n,1,h} \;$
 the characteristic is $\; \chi(x) = 2 \;$ when $\; h \;$ is odd and $\; \chi(x) = 0 \;$ when $\; h \;$ is even,
 i.e. $\; \chi(x) = (-1)^{h+1} + 1 \,$.
 We now proceed by induction on the length $\; m \,$.
 If $\; x \in T_{n,m+1,h} \,$, then $\; x = H^a U^h H^{m+1} D^h H^b \;$
 and hence $\; x = x_1 \vee x_2 \;$ where $\; x_1 = H^a U^h H^m D^h H^{b+1} \,$, $\; x_2 = H^{a+1} U^h H^m D^h H^b \,$.
 Since $\; x_1 \wedge x_2 = H^{a+1} U^h H^{m-1} D^h H^{b+1} \in T_{n,m-1,h} \,$, we have
 $$
  \chi(x) = \chi(x_1) + \chi(x_2) - \chi(x_1\wedge x_2 ) = (-1)^{h+1} m + 1 + (-1)^{h+1} m + 1 - (-1)^{h+1} (m-1) - 1
 $$
 that is $\; \chi(x) = (-1)^{h+1}(m+1)+1 \,$. So, the lemma is proved.
\end{proof}

\begin{remark}
 The maximum $\; \widehat{1} \;$ of $\; \MM_n \;$ is join-irreducible when $\; n = 2 k \;$
 and belongs to $\; T_{2k+1,1,k} \;$ when $\; n = 2 k + 1 \,$.
 Hence, from Lemma \ref{truncpyr}, it follows that $\; \chi(\MM_n) = 1 \;$ when $\; n \;$ is even,
 $\; \chi(\MM_n) = 0 \;$ when $\; n = 4 k + 1 \;$ and $\; \chi(\MM_n) = 2 \;$ when $\; n = 4 k + 3 \,$.
\end{remark}

\medskip

In any Motzkin lattice each quasi-join-irreducible has a particular join-irreducible decomposition,
coming directly from the definition of quasi-join-irreducible element.
However, for our purposes, another kind of decomposition will be useful in representing quasi-join-irreducibles.
\begin{lemma}\label{Lemma-qji-expansion}
 Every quasi-join-irreducible element $\; x \in \MM_n \;$
 can be expressed as $\; x = s_1 \vee \ldots \vee s_k \,$,
 where each $\; s_i \;$ is either join-irreducible or it belongs to $\; T_{n,m,h} \,$,
 and $\; s_i \wedge s_{i+1}\neq \widehat{0} \,$, for every $\; i < k \,$.
\end{lemma}
\begin{proof}
 Write $\; x = t_1 \vee \cdots \vee t_r \;$ as a join of join-irreducibles
 and group together all the consecutive join-irreducibles
 whose join gives rise to a Motzkin path belonging to some $\; T_{n,m,h} \,$.
\end{proof}

Clearly, the decomposition described in Lemma \ref{Lemma-qji-expansion} is not unique.
However, there is a particular way of performing such a decomposition,
which consists of taking truncated pyramids of maximum length, as in the proof.
Such a decomposition will be called the \emph{Motzkin decomposition} of the quasi-join irreducible $\; x \,$.

Now, we are ready to state and proof the fundamental step in the determination (and combinatorial interpretation)
of the characteristic of $\; \MM_n \,$.
Our main proposition will be preceded by a technical lemma.
\begin{lemma}\label{ultimo}
 Let $\; x \in \MM_n \;$ be a quasi-join-irreducible element
 and $\; x = s_1 \vee \cdots \vee s_k \;$ its Motzkin decomposition.
 Then, for every $\; i < k \,$, $\; s_i \wedge s_{i+1} \;$ is join-irreducible.
\end{lemma}
\begin{proof}
 If $\; s_i \;$ and $\; s_{i+1} \;$ are both join-irreducibles,
 the conclusion follows from the definition of Motzkin decomposition.
 If at least one of the two is an element of some $\; T_{n,m,h} \;$,
 then $\; s_i \wedge s_{i+1} \;$ is equal to the meet of two join-irreducibles.
 Indeed, if $\; s_i \;$ is join-irreducible and $\; s_{i+1} \in T_{n,m,h} \,$,
 then $\; s_i \wedge s_{i+1} = s_i \wedge x \,$,
 where $\; s_{i+1} = x \vee x_1 \vee \cdots \vee x_r \;$ is a join-irreducible decomposition of $\; s_{i+1} \;$
 and $\; x \;$ is the join-irreducible having minimum abscissa.
 The remaining cases can be dealt with in a similar way.
\end{proof}

Let $\; x \in \MM_n \,$.
We will write $\; o(x) \;$ for the number of horizontal steps at odd height
and $\; e(x) \;$ for the number of horizontal steps at even nonzero height
(i.e. at even height and not lying on the $x$-axis).
\begin{proposition}
 Let $\; x \in \MM_n \;$ be a quasi-join-irreducible.
 Then $\; \chi(x) = o(x) - e(x) + 1 \,$.
\end{proposition}
\begin{proof}
 Let $\; x = s_1 \vee \cdots \vee s_k \;$ be the Motzkin decomposition of $\; x \,$.
 From formula (\ref{dasilva}) we have
 \begin{equation}\label{motzkindasilva}
  \chi(x) = \chi(s_1 \vee \cdots \vee s_k)
          = \sum_{S\subseteq [k]\atop S\neq \emptyset} (-1)^{|S|-1}\chi \left( \bigwedge_{i\in S}s_i \right) .
 \end{equation}
 Since $\; \MM_n \;$ is a Dyck-like lattice, we can proceed as in Proposition \ref{charqji}.
 We first observe that, if $\; \min S = \max S \,$,
 then the contribution to the sum is $\; \chi(s_1) + \cdots + \chi(s_k) \,$.
 Some of the $\; s_i$'s are join-irreducibles (and so their contribution is $\; 1 \,$),
 but some of them could be paths in $\; T_{n,m,h} \,$.
 If $\; s_{i} \in T_{n,m,h} \,$, then from Lemma \ref{truncpyr}
 it follows that $\; \chi(s_{i}) = (-1)^{h+1}\; m + 1 \,$.
 Starting from this remark, it is not difficult to show that $\; \chi(s_1) + \cdots + \chi(s_k) = k + o(x) - e(x) \,$.
 Now, if $\; \max S = \min S + 1 $, the contribution to the sum is
 $\; -\sum_{i<k} \chi(s_i \wedge s_{i+1}) = -(k - 1) \,$,
 since all meets $\; s_i \wedge s_{i+1} \;$ are necessarily join-irreducible by Lemma \ref{ultimo}.
 Finally, if $\; \max S = \min S + r \,$, with $\; r > 1 \,$,
 using an argument completely analogous to the one used in the Dyck case,
 we find that the contribution of these subsets to the sum in (\ref{motzkindasilva}) is zero.
 So, in conclusion, we have $\; \chi(x) = k + o(x) - e(x) - k + 1 = o(x) - e(x) + 1 \,$.
\end{proof}

Since Motzkin lattices are Dyck-like, from Proposition \ref{qjid} it follows
that every element of a Motzkin lattice has a quasi-join-irreducible decomposition.
Let $\; \| x\| \;$ be the number of all quasi-join-irreducibles in a decomposition of $\; x \in \MM_n \,$,
and let $\; o'(x) \;$ be the number of horizontal steps at odd height different from $\; 1 \,$.
\begin{theorem}\label{motzchar}
 The characteristic of a Motzkin path $\; x \;$ is $\; \chi(x) = \| x \| + o'(x) - e(x) \,$.
\end{theorem}
\begin{proof}
 If $\; x = x_1 \vee \cdots \vee x_k \;$ is a quasi-join-irreducible decomposition of $\; x \,$, then
 $\; \chi(x) = \chi(x_1) + \cdots + \chi(x_k) = o(x_1) - e(x_1) + 1 + \cdots + o(x_k) - e(x_k) + 1 = o'(x) - e(x) + k \,$.
\end{proof}

Also in this case, we have the following remarkable consequence.

\begin{corollary} Two quasi-join-irreducible decompositions of the same
Motzkin path have the same number of elements.
\end{corollary}

A \emph{reverse truncated pyramid} of height $\; h \;$ is
any factor of a Motzkin path of the form $\; D H^k U \,$, where $\; k \geq 1 \;$
and the sequence $\; H^k \;$ of horizontal steps lies on the line $\; y = h \,$.
The height of a peak is given by its ordinate.
If $\; f_h(x) \,$, $\; p_h(x) \,$, $\; t_h(x) \;$ and $\; r_h(x) \;$ are respectively the number
of truncated pyramids, peaks, tunnels and reverse truncated pyramid of height $\; h \;$ in $\; x \,$,
then Theorem \ref{motzchar} can be interpreted combinatorially as stated in
\begin{theorem}
 The characteristic of a Motzkin path $\; x \in \MM_n \;$ is
 $$ \chi(x) = o(x) - e(x) + t_1(x) + f_1(x) + p_1(x) - r_1(x) \, . $$
\end{theorem}
\begin{proof}
 We have only to give a combinatorial interpretation of the term $\; \| x\| \,$.
 The quasi-join-irreducibles in a decomposition of a Motzkin path $\; x \;$ can be of two types only.
 They can be Motzkin paths with a unique elevated factor with no horizontal steps at height $\; 1 \,$,
 and then their number is equal to the number of $1$-tunnels of the given path.
 Otherwise they can be Motzkin paths with a unique peak at height $\; 1 \,$.
 In this case, we have an isolated peak in $\; x \,$,
 or a sequence of consecutive peaks whose abscissas differ by 1.
 In this last case, what we see in $\; x \;$ is a sequence of horizontal steps at height $\; 1 \,$.
 We have three possible configurations:
 \begin{itemize}
  \item[a)] a truncated pyramid at height $\; 1 \,$, $\; UH^m D \,$,
   obtained as the join of $\; m+1 \;$ peaks at height $\; 1 \,$,
  \item[b)] either $\; UH^m U \;$ or $\; DH^m D \,$:
   in both cases, the configuration is obtained as the join of $\; m \;$ peaks at height $\; 1 \,$,
  \item[c)] a reverse truncated pyramid $\; DH^m U \;$ at height $\; 1 \;$
   obtained as the join of $\; m-1 \;$ peaks at height $\; 1 \,$.
 \end{itemize}
 Hence $\; \| x\| = t_1(x) + p_1(x) + h_1(x) + f_1(x) - r_1(x) \,$,
 where $\; h_1(x) \;$ is the number of horizontal steps of $\; x \;$ at height $\; 1 \,$.
\end{proof}

\section{Schr\"oder lattices}

\subsection{Representation}

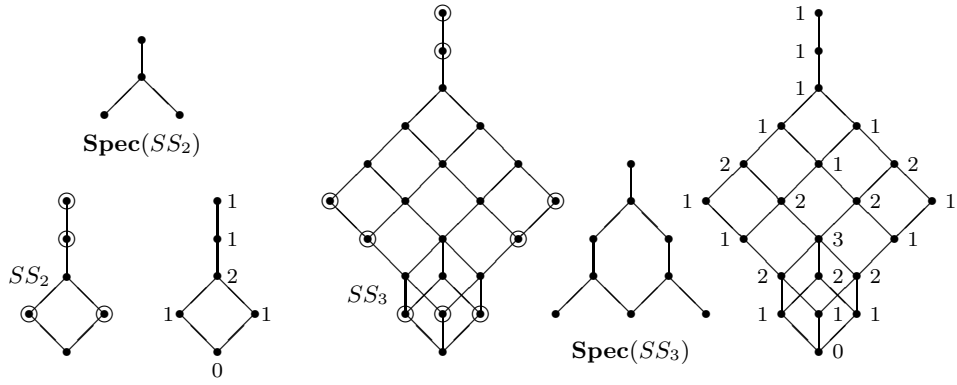
\begin{figure}[h]
\begin{center}
 \setlength{\unitlength}{5mm}
 \begin{picture}(24,9)
  \put(0,0){
   \begin{picture}(2,4)
    \put(1,0){\circle*{0.2}}
    \put(0,1){\circle*{0.2}} \put(0,1){\circle{0.4}}
    \put(2,1){\circle*{0.2}} \put(2,1){\circle{0.4}}
    \put(1,2){\circle*{0.2}}
    \put(1,3){\circle*{0.2}} \put(1,3){\circle{0.4}}
    \put(1,4){\circle*{0.2}} \put(1,4){\circle{0.4}}
    \put(1,0){\line(1,1){1}}
    \put(1,0){\line(-1,1){1}}
    \put(1,2){\line(-1,-1){1}}
    \put(1,2){\line(1,-1){1}}
    \put(1,2){\line(0,1){2}}
    \put(0,2){\makebox(0,0){\small$\SS_2$}}
   \end{picture}}
  \put(4,0){
   \begin{picture}(2,4)
    \put(1,0){\circle*{0.2}} \put(1,-0.5){\makebox(0,0){\footnotesize$0$}}
    \put(0,1){\circle*{0.2}} \put(-0.3,1){\makebox(0,0){\footnotesize$1$}}
    \put(2,1){\circle*{0.2}} \put(2.3,1){\makebox(0,0){\footnotesize$1$}}
    \put(1,2){\circle*{0.2}} \put(1.4,2){\makebox(0,0){\footnotesize$2$}}
    \put(1,3){\circle*{0.2}} \put(1.4,3){\makebox(0,0){\footnotesize$1$}}
    \put(1,4){\circle*{0.2}} \put(1.4,4){\makebox(0,0){\footnotesize$1$}}
    \put(1,0){\line(1,1){1}}
    \put(1,0){\line(-1,1){1}}
    \put(1,2){\line(-1,-1){1}}
    \put(1,2){\line(1,-1){1}}
    \put(1,2){\line(0,1){2}}
   \end{picture}}
  \put(2,6.3){
   \begin{picture}(2,4)
    \put(0,0){\circle*{0.2}}
    \put(2,0){\circle*{0.2}}
    \put(1,1){\circle*{0.2}}
    \put(1,2){\circle*{0.2}}
    \put(0,0){\line(1,1){1}}
    \put(2,0){\line(-1,1){1}}
    \put(1,1){\line(0,1){1}}
    \put(1,-0.8){\makebox(0,0){\small$\Spec(\SS_2)$}}
   \end{picture}}
  \put(8,0){
   \begin{picture}(6,9)
    \put(3,0){\circle*{0.2}}
    \put(2,1){\circle*{0.2}} \put(2,1){\circle{0.4}}
    \put(3,1){\circle*{0.2}} \put(3,1){\circle{0.4}}
    \put(4,1){\circle*{0.2}} \put(4,1){\circle{0.4}}
    \put(2,2){\circle*{0.2}}
    \put(3,2){\circle*{0.2}}
    \put(4,2){\circle*{0.2}}
    \put(1,3){\circle*{0.2}} \put(1,3){\circle{0.4}}
    \put(3,3){\circle*{0.2}}
    \put(5,3){\circle*{0.2}} \put(5,3){\circle{0.4}}
    \put(0,4){\circle*{0.2}} \put(0,4){\circle{0.4}}
    \put(2,4){\circle*{0.2}}
    \put(4,4){\circle*{0.2}}
    \put(6,4){\circle*{0.2}} \put(6,4){\circle{0.4}}
    \put(1,5){\circle*{0.2}}
    \put(3,5){\circle*{0.2}}
    \put(5,5){\circle*{0.2}}
    \put(2,6){\circle*{0.2}}
    \put(4,6){\circle*{0.2}}
    \put(3,7){\circle*{0.2}}
    \put(3,8){\circle*{0.2}} \put(3,8){\circle{0.4}}
    \put(3,9){\circle*{0.2}} \put(3,9){\circle{0.4}}
    \put(3,0){\line(-1,1){1}}
    \put(3,0){\line(0,1){1}}
    \put(3,0){\line(1,1){1}}
    \put(3,1){\line(-1,1){3}}
    \put(3,1){\line(1,1){3}}
    \put(2,2){\line(1,1){3}}
    \put(1,3){\line(1,1){3}}
    \put(0,4){\line(1,1){3}}
    \put(4,2){\line(-1,1){3}}
    \put(5,3){\line(-1,1){3}}
    \put(6,4){\line(-1,1){3}}
    \put(2,1){\line(0,1){1}}
    \put(2,1){\line(1,1){1}}
    \put(4,1){\line(-1,1){1}}
    \put(4,1){\line(0,1){1}}
    \put(3,2){\line(0,1){1}}
    \put(3,7){\line(0,1){2}}
    \put(1,1.5){\makebox(0,0){\small$\SS_3$}}
   \end{picture}}
  \put(14,1){
   \begin{picture}(4,4)
    \put(0,0){\circle*{0.2}}
    \put(2,0){\circle*{0.2}}
    \put(4,0){\circle*{0.2}}
    \put(1,1){\circle*{0.2}}
    \put(3,1){\circle*{0.2}}
    \put(1,2){\circle*{0.2}}
    \put(3,2){\circle*{0.2}}
    \put(2,3){\circle*{0.2}}
    \put(2,4){\circle*{0.2}}
    \put(0,0){\line(1,1){1}}
    \put(2,0){\line(-1,1){1}}
    \put(2,0){\line(1,1){1}}
    \put(4,0){\line(-1,1){1}}
    \put(1,1){\line(0,1){1}}
    \put(3,1){\line(0,1){1}}
    \put(2,3){\line(-1,-1){1}}
    \put(2,3){\line(1,-1){1}}
    \put(2,3){\line(0,1){1}}
    \put(2,-1){\makebox(0,0){\small$\Spec(\SS_3)$}}
   \end{picture}}
  \put(18,0){
   \begin{picture}(6,9)
    \put(3,0){\circle*{0.2}} \put(3.5,0){\makebox(0,0){\footnotesize$0$}}
    \put(2,1){\circle*{0.2}} \put(1.5,1){\makebox(0,0){\footnotesize$1$}}
    \put(3,1){\circle*{0.2}} \put(3.5,1){\makebox(0,0){\footnotesize$1$}}
    \put(4,1){\circle*{0.2}} \put(4.5,1){\makebox(0,0){\footnotesize$1$}}
    \put(2,2){\circle*{0.2}} \put(1.5,2){\makebox(0,0){\footnotesize$2$}}
    \put(3,2){\circle*{0.2}} \put(3.5,2){\makebox(0,0){\footnotesize$2$}}
    \put(4,2){\circle*{0.2}} \put(4.5,2){\makebox(0,0){\footnotesize$2$}}
    \put(1,3){\circle*{0.2}} \put(0.5,3){\makebox(0,0){\footnotesize$1$}}
    \put(3,3){\circle*{0.2}} \put(3.5,3){\makebox(0,0){\footnotesize$3$}}
    \put(5,3){\circle*{0.2}} \put(5.5,3){\makebox(0,0){\footnotesize$1$}}
    \put(0,4){\circle*{0.2}} \put(-0.5,4){\makebox(0,0){\footnotesize$1$}}
    \put(2,4){\circle*{0.2}} \put(2.5,4){\makebox(0,0){\footnotesize$2$}}
    \put(4,4){\circle*{0.2}} \put(4.5,4){\makebox(0,0){\footnotesize$2$}}
    \put(6,4){\circle*{0.2}} \put(6.5,4){\makebox(0,0){\footnotesize$1$}}
    \put(1,5){\circle*{0.2}} \put(0.5,5){\makebox(0,0){\footnotesize$2$}}
    \put(3,5){\circle*{0.2}} \put(3.5,5){\makebox(0,0){\footnotesize$1$}}
    \put(5,5){\circle*{0.2}} \put(5.5,5){\makebox(0,0){\footnotesize$2$}}
    \put(2,6){\circle*{0.2}} \put(1.5,6){\makebox(0,0){\footnotesize$1$}}
    \put(4,6){\circle*{0.2}} \put(4.5,6){\makebox(0,0){\footnotesize$1$}}
    \put(3,7){\circle*{0.2}} \put(2.5,7){\makebox(0,0){\footnotesize$1$}}
    \put(3,8){\circle*{0.2}} \put(2.5,8){\makebox(0,0){\footnotesize$1$}}
    \put(3,9){\circle*{0.2}} \put(2.5,9){\makebox(0,0){\footnotesize$1$}}
    \put(3,0){\line(-1,1){1}}
    \put(3,0){\line(0,1){1}}
    \put(3,0){\line(1,1){1}}
    \put(3,1){\line(-1,1){3}}
    \put(3,1){\line(1,1){3}}
    \put(2,2){\line(1,1){3}}
    \put(1,3){\line(1,1){3}}
    \put(0,4){\line(1,1){3}}
    \put(4,2){\line(-1,1){3}}
    \put(5,3){\line(-1,1){3}}
    \put(6,4){\line(-1,1){3}}
    \put(2,1){\line(0,1){1}}
    \put(2,1){\line(1,1){1}}
    \put(4,1){\line(-1,1){1}}
    \put(4,1){\line(0,1){1}}
    \put(3,2){\line(0,1){1}}
    \put(3,7){\line(0,1){2}}
   \end{picture}}
 \end{picture}
\end{center}
\caption{The Schr\"oder lattices $\; \SS_2 \;$ and $\; \SS_3 \,$, their spectra and the distribution of the characteristic.}
\label{SchroderLatticeS3}
\end{figure}

The Schr\"oder lattice $\; \SS_n \;$ is a distributive lattice of height $\; n^2 \,$,
with minimum $\; H^{2n} \;$ and maximum $\; U^n D^n \,$,
with $\; n - 1 \;$ atoms of the form $\; H^k U D H^{n-k-2} \;$ and one coatom $\; U^{n-1} H^2 D^{n-1} \,$,
The socle is the path $\; s = (UD)^n \;$
and the principal ideal $\; \downarrow s \;$ is isomorphic to a Boolean algebra $\; B_n \,$.
This time, however, the principal filter $\; \uparrow s \;$
just contains an isomorphic copy of the Schr\"oder lattice $\; \SS_{n-1} \;$
consisting of the principal filter generated by the path $\; s' = UH^{n-2}D \,$.
Moreover, the interval $\; [s,s'] \;$ is a Boolean algebra $\; B_{n-2} \;$ (see Figure \ref{SchroderLatticeS3}).
The rank function is given by the area under the path (as for Motzkin lattices).

\begin{proposition}\label{jischr}
 The join-irreducibles of the Schr\"oder lattice $\; \SS_n \;$ are the paths of the form
 \begin{center}
 \setlength{\unitlength}{3mm}
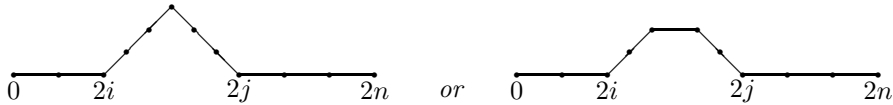
 \begin{picture}(16,4)(0,-1)
  \put(0,0){\circle*{0.2}}
  \put(2,0){\circle*{0.2}}
  \put(4,0){\circle*{0.2}}
  \put(5,1){\circle*{0.2}}
  \put(6,2){\circle*{0.2}}
  \put(7,3){\circle*{0.2}}
  \put(8,2){\circle*{0.2}}
  \put(9,1){\circle*{0.2}}
  \put(10,0){\circle*{0.2}}
  \put(12,0){\circle*{0.2}}
  \put(14,0){\circle*{0.2}}
  \put(16,0){\circle*{0.2}}
  \put(0,0){\line(1,0){4}}
  \put(10,0){\line(1,0){6}}
  \put(4,0){\line(1,1){3}}
  \put(7,3){\line(1,-1){3}}
  \put(0,-0.7){\makebox(0,0){$0$}}
  \put(4,-0.7){\makebox(0,0){$2i$}}
  \put(10,-0.7){\makebox(0,0){$2j$}}
  \put(16,-0.7){\makebox(0,0){$2n$}}
 \end{picture}
 \qquad\textrm{or}\qquad
 \begin{picture}(16,4)(0,-1)
  \put(0,0){\circle*{0.2}}
  \put(2,0){\circle*{0.2}}
  \put(4,0){\circle*{0.2}}
  \put(5,1){\circle*{0.2}}
  \put(6,2){\circle*{0.2}}
  \put(8,2){\circle*{0.2}}
  \put(9,1){\circle*{0.2}}
  \put(10,0){\circle*{0.2}}
  \put(12,0){\circle*{0.2}}
  \put(14,0){\circle*{0.2}}
  \put(16,0){\circle*{0.2}}
  \put(0,0){\line(1,0){4}}
  \put(10,0){\line(1,0){6}}
  \put(4,0){\line(1,1){2}}
  \put(6,2){\line(1,0){2}}
  \put(8,2){\line(1,-1){2}}
  \put(0,-0.7){\makebox(0,0){$0$}}
  \put(4,-0.7){\makebox(0,0){$2i$}}
  \put(10,-0.7){\makebox(0,0){$2j$}}
  \put(16,-0.7){\makebox(0,0){$2n$}}
 \end{picture}
\end{center}
 In particular, $\mathbf{Spec}(\mathcal{S}_n)\simeq \{ (i,j,k)\in \mathcal{C}_n^2 \times \mathcal{C}_1 :j-i\geq 2+2k\}$.
\end{proposition}
\begin{proof}
 Given a Schr\"oder path, there are only two possible
 ways of getting a path which is covered by the starting one:
 either replace an occurrence of $\; U D \;$ with a double horizontal step
 or replace a double horizontal step with $\; D U \,$.
\end{proof}

There are at least two further ways of describing the poset $\; \Spec(\mathcal{S}_n) \,$.
They are essentially equivalent, but the first one is expressed in purely algebraic language
whereas the second one can be considered a sort of combinatorial interpretation.
\begin{enumerate}
 \item
  The \emph{lexicographic product} $\; P \circ Q \;$ of two posets $\; P \;$ and $\; Q \;$
  is the set $\; P \times Q \;$ endowed with the order defined by setting $\; (x_1,y_1) \leq (x_2,y_2) \;$
  when $\; x_1 < x_2 \;$ or $\; x_1 = x_2 \;$ and $\; y_1 \leq y_2 \,$.
  Then $\; \Spec(\mathcal{S}_n) \;$ is isomorphic to the poset
  obtained by $\; \Spec(\DD_{n+1}) \circ \CC_1 \;$ when all the minimal elements are removed
  (see Figure \ref{SchroderLatticeS3}).
 \item
  Denote by $\; \overrightarrow{\Int}(\CC_n ) \;$ the set of oriented intervals of $\; \CC_n \,$.
  An interval $\; I \;$ of a poset $\; P \;$ is said to be \emph{oriented upward (downward)}
  when its elements are listed in such a way that,
  if $\; x < y \;$ in $\; P \,$, then $\; x \;$ precedes (follows) $\; y \;$ in the above listing
  (in this way $\; I \;$ is not just a set but, more precisely, an ordered $t$-uple).
  If $P$ is a chain, $\; P = \CC_n \,$, then its intervals can have only two orientations,
  say \emph{negative} (from top to bottom) and \emph{positive} (from bottom to top).
  We can introduce a partial order on $\; \overrightarrow{\Int}(\CC_n) \;$
  by simply declaring that, for any $\; I, J \in \overrightarrow{\Int}(\CC_n) \,$,
  $\; I \leq J \;$ when either $\; I \;$ is contained in $\; J \;$ as ordinary intervals
  or $\; I = J \;$ as ordinary intervals but $\; I \;$ is negative and $\; J \;$ is positive.
  It is not difficult to see that, endowed with this partial order,
  $\; \overrightarrow{\Int}(\CC_{n-1}) \simeq \Spec(\SS_n) \;$
  (observe that, in $\; \overrightarrow{\Int}(\CC_n ) \,$,
  as far as singleton intervals are concerned,
  there is no distinction between positive and negative intervals).
\end{enumerate}

\begin{proposition}
 The Schr\"oder lattices $\; \mathcal{S}_n \;$ are special Dyck-like lattices where, in particular,
 the meet of any two join-irreducibles is $\; \widehat{0} \;$ or a join-irreducible having a unique peak.
\end{proposition}
\begin{proof}
 As we did for Dyck paths, label each join-irreducible of $\; \SS_n \;$ with the abscissa of its unique (truncated) pyramid:
 this labelling satisfies the condition in the definition of a Dyck-like lattice.
 A join-irreducible path in a Schr\"oder lattice has either a unique peak
 or a unique truncated pyramid with a double horizontal step.
 Since two Schr\"oder paths of the same length can cross only at points with integer coordinates
 and cannot cross in the middle point of a double horizontal step of one of them,
 if the intersection of two join-irreducibles is not the minimum $\; \widehat{0} \,$,
 then it is necessarily a join-irreducible with a unique peak.
\end{proof}

\subsection{Characteristic}

The characteristic of Schr\"oder lattices admits a combinatorial interpretation
analogous to the one given for Dyck lattices.
Indeed, following the same lines of the proofs of Proposition \ref{dyckqji} and Corollary \ref{chardyck}
and using Proposition \ref{jischr}, we have
\begin{theorem}
 A Schr\"oder path is quasi-join-irreducible if and only if it has exactly one $0$-tunnel.
\end{theorem}
\begin{theorem}
 The characteristic of a Schr\"oder path equals the number of its $0$-tunnels.
\end{theorem}

\section{Final remarks on rank unimodality}

An interesting property common to several sequences arising in combinatorics is unimodality.
Specifically, a \emph{sequence} $\; \{ a_0, a_1, \ldots, a_n \} \;$ of positive integers is \emph{unimodal}
when there exists an index $\; k \;$ such that
$\; a_0 \leq a_1 \leq \cdots \leq a_k \geq a_{ k + 1 } \geq \cdots \geq a_n \,$,
and a \emph{polynomial} is \emph{unimodal} when the sequence of its coefficients is unimodal.
In the case of ranked posets this is often a property of the distribution of the elements of given rank.
More precisely, a (finite) \emph{poset} $\; P \;$ is \emph{rank unimodal} when it is ranked and
its rank polynomial (or equivalently the sequence of its Whitney numbers) is unimodal.
The \emph{Whitney number} $\; W_k(P) \;$ is the number of all elements of $\; P \;$ with rank $\; k \;$
whereas the \emph{rank polynomial} is $\; R(P;q) = W_0(P) + W_1(P) q + \cdots + W_h(P) q^h \,$,
where $\; h = r(P) \;$ is the height of $\; P \,$.

\paragraph{Dyck lattices.}
Let $\; D_n(q) \;$ be the rank polynomial of $\; \DD_n \,$.
Since every non-empty Dyck path decomposes as $\; U \gamma' D \gamma'' \;$
(where $\; \gamma' \;$ and $\; \gamma'' \;$ are arbitrary Dyck paths), we have the recurrence
\begin{displaymath}
 D_{n+1}(q) = \sum_{k=0}^n q^k D_k(q) D_{n-k}(q)
\end{displaymath}
with the initial condition $\; D_0(q) = 1 \,$.
It is easy to see that the generating series for these polynomials satisfies the identity
$\; D(q,t) = (1 - t D(q,qt))^{-1} \;$ 
from which it is possible to obtain an expansion as a continued fraction \cite{Stanley2}.
The polynomials $\; D_n(q) \;$ define a $q$-analog of Catalan numbers,
namely $\; q^{{n\choose 2}} D_n(1/q) = C_n(q) \,$,
where the $\; C_n(q)$'s are the \emph{$q$-Catalan numbers} defined as the sum $\; \sum_x q^{\AA(x)} \;$
over all lattice paths from $\; (0,0) \;$ to $\; (n,n) \,$,
with steps $\; (1,0) \;$ and $\; (0,1) \,$, never rising above the line $\; y = x \,$,
where $\; \AA(x) \;$ is the area of the region determined by the path and the $x$-axis \cite{FH} \cite[p. 235]{Stanley2}.
The Whitney numbers of $\; \DD_n \;$ appears in \cite{sloane} (essentially) as sequence A129182.
It is still an open problem \cite{BSS} to prove or disprove
that the rank polynomials $\; D_n(q) \;$ are unimodal for every $\; n \in \NN \,$.
This problem is also mentioned in \cite{S}, where it is conjectured that
the Young's lattices associated with the staircase partitions $\; (n,n-1,\ldots ,2,1) \;$ are rank unimodal.

\paragraph{Motzkin lattices.}
Let $\; M_n(q) \;$ be the rank polynomial of $\; \MM_n \,$.
Since every non-empty Motzkin path decomposes as $\; H \gamma \;$ or $\; U \gamma' D \gamma'' \;$
(where $\; \gamma \,$, $\; \gamma' \;$ and $\; \gamma'' \;$ are arbitrary Motzkin paths),
we have the recurrence
\begin{displaymath}
 M_{n+2}(q) = M_{n+1}(q) + \sum_{k=0}^n q^{k+1}M_k(q)M_{n-k}(q)
\end{displaymath}
with the initial conditions $\; M_0(q) = M_1(q) = 1 \,$.
Their generating series $\; M(q,t) \;$ satisfies the identity
$\; M(q,t) = (1-t-qt^2 M(q,qt))^{-1} \;$ 
and hence admits an expansion as a continued fraction.
The Whitney numbers of $\; \MM_n \;$ appear in \cite{sloane} as sequence A129181
and also in this case we can conjecture that the lattices $\; \MM_n \;$ are rank-unimodal.

\paragraph{Schr\"oder lattices.}
Let $\; S_n(q) \;$ be the rank polynomial of $\; \SS_n \,$.
Since every non-empty Schr\"oder path decomposes as $\; HH \gamma \;$ or $\; U \gamma' D \gamma'' \;$
(where $\; \gamma \,$, $\; \gamma' \;$ and $\; \gamma'' \;$ are arbitrary Schr\"oder paths), we have the recurrence
\begin{displaymath}
 S_{n+1}(q) = S_n(q) + \sum_{k=0}^n q^{2k+1}S_k(q)S_{n-k}(q)
\end{displaymath}
with the initial condition $\; S_0(q) = 1 \,$.
Their generating series $\; S(q;x) \;$ satisfies the identity
$\; S(q;x) = (1-x-qxS(q;q^2x))^{-1} \;$ 
and hence also this time it has an expansion as a continued fraction.
The Whitney numbers of $\; \SS_n \;$ appear in \cite{sloane} as sequence A129179.
Also in this case it is still an open problem \cite{BSS} to prove the rank-unimodality of the lattices $\; \SS_n \,$.

\end{document}